\documentclass[a4paper,12pt]{article}

\usepackage[utf8]{inputenc}
\usepackage[T1]{fontenc}

\usepackage{amsmath,
            mleftright,
            comment,
            amssymb,
            amsthm,
            nicefrac,
            xcolor,
            enumerate,
            authblk,
            color,
            etoolbox,
            mathrsfs,
            mathabx,
            bbm,
            xparse,
            stmaryrd,
            geometry,
            enumitem,
            mathtools
            }

\mleftright   

\usepackage[sort,nocompress]{cite}

\usepackage[english]{datetime2}
            
\usepackage[colorlinks=true]{hyperref}
            
\newcommand{\N}{\ensuremath{\mathbb{N}}}

\newcommand{\R}{\ensuremath{\mathbb{R}}}
\newcommand{\C}{\ensuremath{\mathbb{C}}}

\newcommand{\1}{\mathbbm{1}}
\newcommand{\dx}{\, {\rm d}}
\newcommand{\logNorm}{\mu}
\newcommand{\idMatrix}{\operatorname{I}}

\newcommand{\dpp}{\text{.}}
\newcommand{\dc}{\text{,}}
\newcommand{\Ts}{T_{\text{s}}}

\newcommand{\U}{{\bf U}}
\newcommand{\LL}{\mathbb{L}}
\newcommand{\V}{\mathbb{V}}
\newcommand{\CC}{\mathcal{C}}
\newcommand{\semi}{\mathbb{T}}
\newcommand{\temp}[1]{v^{(#1)}}
\newcommand{\0}{{\bf 0}}
\newcommand{\fh}{\mathfrak{h}}

\newcommand{\fc}{\mathfrak{c}}
\newcommand{\fx}{\mathfrak{x}}
\newcommand{\w}[1]{\mathbb{W}^{(#1)}}
\newcommand{\mm}{\mathbb{M}}

\usepackage[colorinlistoftodos]{todonotes}
\tikzset{notestyleraw/.append style={align=justify}}
                        
\usepackage[english]{datetime2}

\newcounter{todocounter}

\newcommand{\smallsum}{\textstyle\sum}
\newcommand{\SmallSum}{\textstyle\sum\limits}

\newcommand{\induct}{\dashrightarrow}
\newcommand{\with}{\curvearrowleft}
\newcommand{\lrSpace}{\ensuremath{\mkern-1.5mu}}

\DeclarePairedDelimiter{\pr}{(}{)}
\DeclarePairedDelimiter{\br}{[}{]}
\DeclarePairedDelimiter{\cu}{\{}{\}}
\DeclarePairedDelimiter{\abs}{\lvert}{\rvert}
\DeclarePairedDelimiter{\norm}{\lVert}{\rVert}

\DeclarePairedDelimiter{\vt}{\langle}{\rangle}

\allowdisplaybreaks

\usepackage[sort,capitalize]{cleveref}

\crefformat{equation}{(#2#1#3)}
\crefname{equation}{}{}
\crefname{enumi}{item}{items}
\crefname{subsection}{Subsection}{Subsections}
\crefname{setting}{Setting}{Settings}

\newtheorem{theorem}{Theorem}[section]

\theoremstyle{definition}
\newtheorem{definition}[theorem]{Definition}
\newtheorem{setting}[theorem]{Setting}
\newtheorem{remark}[theorem]{Remark}

\numberwithin{equation}{section}

\ExplSyntaxOn

\NewDocumentCommand{\enum}{ O{;} m o }
 {
  \my_enum:nnn { #1 } { #2 } { #3 }
 }

\seq_new:N \l__my_enum_seq
\tl_new:N \l__my_enum_item_tl
\prop_const_from_keyval:Nn \l__verbs
{
show = shows ,
imply = implies , 
demonstrate = demonstrates ,
prove = proves ,
establish = establishes ,
ensure = ensures ,
assure = assures
}
\int_new:N \l__number_of_args

\cs_new_protected:Nn \my_enum:nnn
 {
  \seq_set_split:Nnn \l__my_enum_seq { #1 } { #2 }
  \seq_remove_all:Nn \l__my_enum_seq {}
  \int_set_eq:NN \l__number_of_args { \seq_count:N \l__my_enum_seq }
  \seq_use:Nnnn \l__my_enum_seq { ~and~ } { ,~ } { ,~and~ }
  \IfNoValueTF{#3}{}{
    \space
    \int_compare:nNnTF{ \l__number_of_args } < {2}{ \prop_item:Nn \l__verbs {#3} }{ #3 }
  }
 }

\seq_new:N \g_cflist_loaded
\seq_new:N \g_cflist_pending

\NewDocumentCommand{\cfadd}{ m }
{
  \seq_if_in:NnF \g_cflist_loaded { #1 } {
    \seq_if_in:NnF \g_cflist_pending { #1 } {
      \seq_gput_right:Nn \g_cflist_pending { #1 }
    }
  }
}

\NewDocumentCommand{\cfconsiderloaded}{ m }{
  \seq_gput_right:Nn \g_cflist_loaded {#1}
}

\NewDocumentCommand{\cfremove}{ m }
{
  \seq_gremove_all:Nn \g_cflist_pending { #1 }
}

\NewDocumentCommand{\cfload}{ o }
{
  \seq_if_empty:NTF \g_cflist_pending {\unskip} {
    (cf.\ \cref{\seq_use:Nn \g_cflist_pending {,}})\IfValueTF{#1}{#1~}{\unskip}
    \seq_gconcat:NNN \g_cflist_loaded \g_cflist_loaded \g_cflist_pending
    \seq_gclear:N \g_cflist_pending
  }
}

\NewDocumentCommand{\cfclear} {} {
  \seq_gclear:N \g_cflist_loaded
  \seq_gclear:N \g_cflist_pending
}

\NewDocumentCommand{\cfout}{ o }
{
  \seq_if_empty:NTF \g_cflist_pending {\unskip} {
    (cf.\ \cref{\seq_use:Nn \g_cflist_pending {,}})\IfValueTF{#1}{#1~}{\unskip}
    \seq_gclear:N \g_cflist_pending
  }
}

\NewDocumentCommand{\ifnocf} { m } {
  \seq_if_empty:NT \g_cflist_pending { #1 }
}

\NewDocumentEnvironment {athm} {m m} {%
\begin{#1}\label{#2}\global\def\loc{#2}%
}{%
\end{#1}%
}

\NewDocumentEnvironment{aproof} {} {%
\begin{proof}[Proof~of~\cref{\loc}]%
}{%
\finishproofthus
\end{proof}%
}

\ExplSyntaxOff

\newcommand{\finishproofthus}{The proof of \cref{\loc} is thus complete.}

\ExplSyntaxOn

\NewDocumentEnvironment{flexmath}{ m o }{
  \str_if_eq:noTF {a} {#1} {
    \begin{equation}
    \IfValueT{#2}{\label{eq:\loc.#2}}
    \begin{aligned}
  } {
    \catcode`&=9
    \renewcommand{\\}{}
    \str_if_eq:noTF {d} {#1} {
      \begin{equation}
      \IfValueT{#2}{\label{eq:\loc.#2}}
    } {
      \begin{math}
    }
  }
}{
  \str_if_eq:noTF {i} {#1} {
    \end{math}
    \catcode`&=4
  } {
    \str_if_eq:noTF {d} {#1} {
      \end{equation}
    } {
      \end{aligned}
      \end{equation}
    }
  }
}

\ExplSyntaxOff

\begin{document}

\title{A positivity- and monotonicity-preserving \\ 
nonlinear operator splitting approach for \\
approximating solutions to quenching-combustion semilinear partial differential equations}

\author{
Joshua Lee Padgett$^{1,2}$ and Eduardo Servin$^{3}$
\bigskip
\\
\small{$^1$ Department of Mathematical Sciences, University of Arkansas,}
\vspace{-0.1cm}\\
\small{Fayetteville, Arkansas 72701, USA, e-mail: \texttt{padgett@uark.edu}}
\smallskip
\\
\small{$^2$ Center for Astrophysics, Space Physics, and Engineering Research,}
\vspace{-0.1cm}\\
\small{Baylor University, Waco, Texas 76798, USA, e-mail: \texttt{padgett@uark.edu}}
\smallskip
\\
\small{$^3$ Department of Mathematics, Baylor University,}
\vspace{-0.1cm}\\
\small{Waco, Texas 76798, USA, e-mail: \texttt{eduardo\_servin1@baylor.edu}}
\smallskip
}

\date{\today}

\maketitle

\begin{abstract}
In recent years, there has been a large increase in interest in numerical algorithms which preserve various qualitative features of the original continuous problem.
Herein, we propose and investigate a numerical algorithm which preserves qualitative features of so-called quenching combustion partial differential equations (PDEs).
Such PDEs are often used to model solid-fuel ignition processes or enzymatic chemical reactions and are characterized by their singular nonlinear reaction terms and the exhibited positivity and monotonicity of their solutions on their time intervals of existence.
In this article, we propose an implicit nonlinear operator splitting algorithm which allows for the natural preservation of these features.
The positivity and monotonicity of the algorithm is rigorously proven.
Furthermore, the convergence analysis of the algorithm is carried out and the explicit dependence on the singularity is quantified in a nonlinear setting.
\end{abstract}

\tableofcontents

\section{Introduction}

So-called quenching-combustion partial differential equations (PDEs) arise in the modeling of highly sophisticated, yet important, natural phenomena where singularities may develop as the solution evolves in time.
Such singularities often result from the energy of a system concentrating and approaching its activation criterion \cite{Bebernes_89,Poin,MR2604963,MR1955412}. 
These quenching-combustion PDEs have been studied extensively and are often referred to as \emph{Kawarada PDEs}, after Hideo Kawarada who originally pioneered their research \cite{Kawa}.

Quenching-combustion PDEs are often best understood via the solid-fuel ignition problems they are used to model; see, for instance, \cite{MR1343330}. 
Consider a typical solid-fuel ignition process which occurs in an idealized combustion chamber. 
If the combustion chamber is filled with fuel and air which are appropriately mixed, then prior to ignition the temperature in the chamber may increase monotonically until a certain critical value is reached. 
However, the rates of such temperature changes can occur in a nonlinear manner throughout the media. 
This nonuniform distribution of heat may result in high temperatures being extremely localized within the combustion chamber and ultimately lead to an ignition once the peak temperature reaches a certain threshold (cf., e.g., \cite{Bebernes_89,Poin}). 
This phenomenon is carefully characterized by quenching-combustion PDE models in which the temporal derivative of the solution may grow at an explosive rate, while the solution itself remains bounded. 
This strong nonlinear singularity, which is referred to as a \emph{quenching singularity}, is a feature that any relevant model must capture. 

Theoretical properties of quenching-combustion problems are well-understood (cf., e.g., \cite{MR1635771,MR990863,MR1042837,MR969520,MR3852609} and the references therein).
While some quenching-combustion PDEs can be viewed as being related to the well-known class of \emph{blow-up} PDEs (cf., e.g., 
\cite{MR2422113,MR1374282,MR2084210,MR1897690} and the references therein)
it is the case that quenching-combustion PDEs exhibit novel computational difficulties due to the observed blow-up of the temporal derivative while the solution itself remains bounded
(cf., e.g., \cite{Sheng5,Sheng8,Chen,Josh2,Josh1,Josh3,Josh_Thesis,Chan2} and the references therein).
The aforementioned numerical approximation approaches have been shown to provide acceptable numerical results (e.g., they have demonstrated appropriate accuracy outside of some neighborhood of the singularity), but the theoretical analyses of the proposed methods are still lacking.
In particular, the standard approach when studying numerical approximations of quenching-combustion PDEs is to ``freeze'' the nonlinear term during the classical von Neumann stability analysis---a technique which significantly limits the impact that the singularity can have on any ensuing analysis.
These results provide insight into the numerical methods prior to quenching, but leave many open questions regarding the effects that quenching may have on the numerical method.
Moreover, it is the case that numerical methods for quenching-combustion PDEs are often designed without the preservation of qualitative features in mind and then simply place overly-restrictive conditions on the underlying spatial grids to recover said features.

It is precisely the subject of this article to construct and analyze a nonlinear operator splitting method to approximate certain quenching-combustion semilinear PDEs and provide the rigorous analysis which demonstrates the positivity, monotonicity, and convergence for the implicit numerical scheme.
To better illustrate the findings of this work, we present the following result, \cref{th:main} below.

\begin{theorem}\label{th:main}
Let $a, T \in (0,\infty)$, $\delta \in (0,\min\{1,T\})$, $N \in \N = \{1,2,3,\ldots\}$, let $h = \nicefrac{2a}{(N+1)}$,
let $u \colon [0,T) \times [-a,a] \to \R$ satisfy for all $t \in [0,T)$, $x \in [-a,a]$ that $u(t,-a) = u(t,a) = u(0,x) = 0$, $\lim_{s\to T^-} [ \sup_{y \in [-a,a]} u(s,y) ] = 1$, and
\begin{equation}\label{pde1}
    \pr[\big]{ \tfrac{\partial}{\partial t} u}\lrSpace(t,x) = \pr[\big]{ \tfrac{\partial^2}{\partial x^2} u}\lrSpace(t,x) + \pr[\big]{ 1-u(t,x) }^{-1} \dc
\end{equation}
let 
$\tau_{-1},\tau_0,\tau_1,\tau_2,\ldots \in [0,T)$,
$\temp{0},\temp{1},\temp{2},\ldots \in \R^N$ satisfy for all $k \in \N_0 = \N \cup \{0\}$, $n \in \{1,2,\ldots,N\}$ that $\tau_{-1}=1$, $(\temp{0})_n = 0$, $\tau_k 
= \delta \br[]{ \textstyle\min_{(i,j)\in\{1,2,\ldots,N\}\times\{k,k+1\}} \pr[]{ 1 - \temp{j+1}_i }^2 }
$, and
\begin{equation}\label{split1}
\begin{split}
    & (\temp{k+1})_n 
    - (\temp{k})_n - \tau_k\pr[\big]{ 1-(\temp{k+1})_n }^{-1}
    \\
    & \quad = 
    \tfrac{\tau_k}{h^2} \br[\Big]{ \1_{(1,N]}(n) (\temp{k+1})_{\max\{0,n-1\}} -2(\temp{k+1})_n + \1_{[1,N)}(n) (\temp{k+1})_{\min\{N,n+1\}} } \\
    & \qquad - \tfrac{(\tau_k)^2}{h^2} \Bigl[ \1_{(1,N]}(n) \pr[\big]{1-(\temp{k+1})_{\max\{0,n-1\}}}^{-1} -2\pr[\big]{1-(\temp{k+1})_n}^{-1} \\
    & \qquad \qquad + \1_{[1,N)}(n) \pr[\big]{1-(\temp{k+1})_{\min\{N,n+1\}}}^{-1} \Bigr] \dpp
\end{split}
\end{equation}
Then
\begin{enumerate}[label=(\roman*)]
    \item\label{th:main_1} it holds for all $t_1,t_2 \in [0,T)$, $x \in [-a,a]$ with $t_1 \le t_2$ that $0 \le u(t_1,x) \le u(t_2,x) < 1$,
    \item\label{th:main_4} it holds for all $k \in \N_0$, $n \in \{1,2,\ldots,N\}$ that $0 \le (\temp{k})_n \le (\temp{k+1})_n < 1$, 
and
    \item\label{th:main_6} there exist $\mathfrak{C}_0, \mathfrak{C}_1, \mathfrak{C}_2, \ldots \in \R$ such that for all $k \in \N_0$ with $\sum_{j=0}^k \tau_j < T$ it holds that
    \begin{equation}
\pr[\Big]{ \textstyle\sum_{n=1}^N h \abs[\big]{ (\temp{k})_n - u\pr[\big]{ {\textstyle\sum_{j=0}^k \tau_j} , -a+nh } }^2 }^{\!\nicefrac{1}{2}}
\le \mathfrak{C}_k \pr[\big]{ h^2 + \tau_0 } \dpp
    \end{equation}
\end{enumerate}
\end{theorem}

\cref{th:main} is an immediate consequence of \cref{cor:final} in \cref{sec:5} below.
\cref{cor:final}, in turn, follows from \cref{th:final} which is the main result of the article (see \cref{sec:5} below for details).
In the following, we provide some explanatory comments concerning the mathematical objects appearing in \cref{th:main} above.

In \cref{th:main} we intend to approximate the solution of the PDE in \cref{pde1}.
The strictly positive real number $T \in (0,\infty)$ in \cref{th:main} describes the time horizon for the local interval of existence of the PDE in \cref{pde1},
the strictly positive real number $a \in (0,\infty)$ describes the size of the symmetric spatial domain,
and the function $u\colon [0,T) \times [-a,a] \to \R$ is the solution of the PDE in \cref{pde1}.
The condition that $\lim_{s\to T^-} [ \sup_{y \in [-a,a]} u(s,y) ] = 1$ ensures that $a$ is large enough to guarantee that the solution to the PDE in \cref{pde1} only exists locally.
Moreover, these conditions guarantee that the solution to the PDE in \cref{pde1} is both positive and monotonically increasing, which is summarized in \cref{th:main_1} in \cref{th:main}.

The natural number $N \in \N$ in \cref{th:main} describes the number of internal grid points employed to discretize the interior of the domain $[-a,a]$ via a standard second-order finite difference approximation.
The positive real number $h \in (0,\infty)$ in \cref{th:main} describes the width of the aforementioned intervals.
For simplicity, we have used a uniform spatial grid in \cref{th:main} (which is a special case of the nonuniform grids used in \cref{th:final,cor:final}).
The real numbers $\tau_0,\tau_1,\tau_2,\ldots \in [0,T)$ in \cref{th:main} describe the set of non-uniform temporal step sizes used to iterate our proposed numerical method forward in time, with $\tau_{-1} =1$ simply serving as a convenient initializing variable.
The vectors $\temp{0},\temp{1},\temp{2},\ldots \in \R^N$ in \cref{th:main} represent the numerical approximation of the solution to \cref{pde1} obtained via the nonlinear splitting algorithm in \cref{split1}.
The condition in \cref{th:main} that $\delta,\tau_{-1},\tau_0,\tau_1,\tau_2,\ldots \in [0,T)$ satisfy for all $k \in \N_0$ that $\tau_{-1} = 1$ and $\tau_k 
= \delta \br[]{ \textstyle\min_{(i,j)\in\{1,2,\ldots,N\}\times\{k,k+1\}} \pr[]{ 1 - \temp{j+1}_i }^2 }
$, where $\delta \in (0,1)$ is some chosen tolerance, guarantees that the numerical approximations in \cref{split1} preserve the positivity and monotonicity of the solution to the PDE in \cref{pde1}.
These aforementioned conditions are outlined in \cref{th:main_4} in \cref{th:main}.
\Cref{th:main_6} in \cref{th:main} shows that the numerical approximation in \cref{split1} is a numerical method which exhibits a first-order convergence rate in time, up to the quenching singularity.

The remainder of this article is structured as follows.
In \cref{sec:2} we introduce the general quenching-combustion semilinear PDEs and demonstrate some of their theoretical properties.
Afterwards, in \cref{sec:3}, we introduce a semidiscretized system which approximates the continuous quenching-combustion PDEs. 
The positivity, monotonicity, and convergence of the approximation is rigorously demonstrated.
Next, in \cref{sec:4}, we introduce the proposed nonlinear operator splitting method for approximating general quenching-combustion PDEs and demonstrate the desired positivity and monotonicity.
Finally, in \cref{sec:5}, we prove the necessary convergence results for the proposed nonlinear splitting algorithm in a fully nonlinear setting.

\section{Quenching-combustion PDEs}\label{sec:2}

In this section we first introduce the general class of quenching-combustion PDEs which will be of interest throughout this article in \cref{setting1} below.
Afterwards, we provide some clarifying remarks on the abstract mathematical objects presented in \cref{setting1} below.
In particular, we provide an example of a nonlinear function $f \colon [0,1) \to \R$ which would satisfy the necessary properties for quenching to occur (and occurs quite frequently in practice).
Next, in \cref{lem:quench_props} below, we prove results regarding the regularity, positivity, and monotonicity in time of the solutions of quenching-combustion PDEs.
Finally, we establish in \cref{lem:quench} below the fact that solutions to quenching-combustion PDEs will quench under the assumptions outlined in \cref{setting1} below (see \cref{def:quench} below).

\subsection{Setting}\label{sec:quench_setting}

In this subsection we present the assumptions needed to study a one-dimensional quench\-ing-combustion problem (see \cref{pde2} below).
Generally speaking, we assume that the nonlinearity $f\colon [0,1) \to \R$ is positive, differentiable, and convex (see \cref{f_cond} below).
Moreover, we assume that $f \colon [0,1) \to \R$ satisfies certain Lipschitz and growth conditions (see \cref{f_cond,f_cond2} below).
The assumption that the initial condition satisfies \cref{pos_init} ensures that the solution to \cref{pde2} remains positive and increases monotonically on its time interval of existence.

\begin{setting}\label{setting1}
Let $a,T \in (0,\infty)$, 
let $L \colon [0,1)\times[0,1) \to [0,\infty)$, $f \in C^1([0,1),\R)$ satisfy for all $x , y\in [0,1)$, $s \in [0,1]$
that $f(0) > 0$, $f'(x) > 0$,
\begin{equation}\label{f_cond}
\abs{f(x)-f(y)} \le L_{x,y} \abs{x-y} \dc
\qquad
f\pr[\big]{sx + (1-s)y} \le sf(x) + (1-s)f(y) \dc
\end{equation}
\begin{equation}\label{f_cond2}
\lim_{w\to 1^-} f(w) = \infty \dc
\quad \text{and} \quad
\int_0^1 f(w) \dx w = \infty \dc
\end{equation}
let $u_0 \in C^2( [-a,a] , [0,1) )$ satisfy for all $x \in [-a,a]$ that $u_0(-a) = u_0(a) = 0$ and
\begin{equation}\label{pos_init}
\pr[\big]{ \tfrac{\partial^2}{\partial x^2} u_0}\lrSpace(x) + f(u_0(x)) > 0 \dc
\end{equation}
and let $ u \colon [0,T) \times [-a,a] \to \R$ satisfy for all $t \in [0,T)$, $x \in [-a,a]$ that
$u(t,-a) = u(t,a) = 0$, $u(0,x) = u_0(x)$, 
$\lim_{s\to T^-} [ \sup_{y\in[-a,a]} u(s,y) ] = 1$, 
and
\begin{equation}\label{pde2}
\pr[\big]{ \tfrac{\partial}{\partial t} u}\lrSpace(t,x) = \pr[\big]{ \tfrac{\partial^2}{\partial x^2} u}\lrSpace(t,x) + f(u(t,x)) \dpp
\end{equation}
\end{setting}

\begin{remark}
Note that \cref{f_cond} in \cref{setting1} guarantees that $f$ is \emph{locally Lipschitz} (cf., e.g., Marsden \cite[page 161, Theorem 3]{MR0357693}). In particular, one has that for all $x,y\in[0,1)$ it holds that
\begin{equation}
L_{x,y} \le \sup\left( \left\{ \frac{\abs{f(w)-f(z)}}{\abs{w-z}} \colon w,z \in [0,\max\{x,y\}) \right\} \cup \{0\} \right) < \infty \dpp
\end{equation}
\end{remark}

\begin{remark}
Observe that the assumption in \cref{setting1} that
\begin{equation}\label{quench_cond_a}
\textstyle\lim_{s\to T^-} \left[ \sup_{y\in[-a,a]} u(s,y) \right] = 1
\end{equation}
is a condition on $T$, $a$, and $u$. 
It is often the case that authors will choose a value $a \in (0,\infty)$ and then the goal of their study is to deduce whether or not there exists $T_a \in (0,\infty)$ for which $u$ satisfies \cref{quench_cond_a} (cf., e.g., \cite{Josh3,MR2363217,MR586912,MR2045056} and the references therein).
Since the determination of such relationships is not the goal of the current study, \cref{setting1} simply assumes that the arbitrary values $a$ and $T$ are chosen to guarantee that the solution $u$ to \cref{pde2} satisfies \cref{quench_cond_a}.
\end{remark}

We close \cref{sec:quench_setting} by briefly outlining some details regarding a quintessential example of a quenching-combustion PDE.
Let $f \colon [0,1) \to \R$ satisfy for all $x\in[0,1)$ that $f(x) = (1-x)^{-1}$.
Note that for all $x,y\in[0,1)$ it holds that $f(0) = 1 > 0$, $f'(x) = (1-x)^{-2} > 0$, and
\begin{equation}
\abs{f(x) - f(y)} 
= \abs[\big]{ (1-x)^{-1} - (1-y)^{-1} } 
\le \left[ \frac{1}{\pr[\big]{ 1-\max\{x,y\} }^2 } \right] \abs{x-y} \dpp
\end{equation}
In addition, observe that the fact that $f$ is a convex function ensures that for all $x,y\in[0,1)$, $s\in[0,1]$ it holds that
\begin{equation}
f\pr[\big]{ sx + (1-s)y } \le sf(x) + (1-s)f(y) \dpp
\end{equation}
Moreover, note that
\begin{equation}
\lim_{w\to 1^-} f(w) = \lim_{w\to 1^-} (1-w)^{-1} = \infty
\end{equation}
and
\begin{equation}
\int_0^1 f(w) \dx w = \int_0^1 (1-w)^{-1} \dx w 
= \lim_{v \to 1^-} \pr[\big]{ -\ln\abs{1-v} } = \infty \dpp
\end{equation}
Hence, $f$ satisfies the conditions outlined in \cref{setting1} (i.e., \cref{f_cond,f_cond2}).
Next, let $u_0 \colon [-\sqrt{2},\sqrt{2}] \to [0,1)$ satisfy for all $x\in[-\sqrt{2},\sqrt{2}]$ that $u_0(x) = 0$. 
Observe that $u_0$ and $f$ clearly satisfy \cref{pos_init}.
Then, e.g., Kawarada \cite[Theorem]{Kawa} and, e.g., Kawarada \cite[Lemma]{Kawa} assure that
there exists $T \in (0,\infty)$ such that for all $t\in[0,T)$, $x\in[-\sqrt{2},\sqrt{2}]$ it holds that $u(t,-\sqrt{2}) = u(t,\sqrt{2}) = u(0,x) = 0$, 
\begin{equation}\label{ex_quench_prop}
\textstyle \lim_{s\to T^-} \br[\big]{ \sup_{y \in [-\sqrt{2},\sqrt{2}]} u(x,y) } = \lim_{s\to T^-} u(s,0) = 1 \dc
\end{equation}
and
\begin{equation}
\pr[\big]{ \tfrac{\partial}{\partial t} u}\lrSpace(t,x) = \pr[\big]{ \tfrac{\partial^2}{\partial x^2} u}\lrSpace(t,x) + \pr[\big]{ 1 - u(t,x) }^{-1} \dpp
\end{equation}
Furthermore, it is the case that in this one-dimensional setting, it holds that $T \in (0,\infty)$ depends \emph{only} on the choice of $a$ (in the above, based on \cref{setting1}, it is the case that $a=\sqrt{2}$). 
In fact, it was shown in Chan and Chen \cite[Section 4]{Chen}
that $a \in (0,\infty)$ must satisfy $a \ge 0.7651$ in order to guarantee that \cref{ex_quench_prop} holds.
Interested readers may refer to Padgett and Sheng \cite[Figure 7]{Josh3} for more details on the dynamic relationship between the value of $a \in (0,\infty)$ and $T\in(0,\infty)$ which ensure that \cref{ex_quench_prop} holds.

\subsection{Properties of solutions to quenching-combustion PDEs}\label{sec:quench_props}

In this subsection, we provide theoretical analysis of the solution $u$ to \cref{pde2} under the assumptions outlined in \cref{setting1} above.
The results below exist in the literature in various forms in the case that for all $x\in[-a,a]$ it holds that $u_0(x) = 0$.
Thus, we include the following results for completeness and for the reader's convenience.

\begin{athm}{lemma}{lem:quench_props}
Assume \cref{setting1}.
Then
\begin{enumerate}[label=(\roman*)]
\item\label{lem:quench_props_i0} it holds that $u$ is a unique solution to \cref{pde2},
\item\label{lem:quench_props_i1} it holds that $u \in C^{1,2}([0,T)\times[-a,a],\R)$,
\item\label{lem:quench_props_i2} it holds for all $t \in [0,T)$, $x\in[-a,a]$ that $u(t,x) \in [0,\infty)$, and
\item\label{lem:quench_props_i3} it holds for all $t_1,t_2 \in [0,T)$, $x \in [-a,a]$ with $t_1 \le t_2$ that $u(t_1,x) \le u(t_2,x)$
\end{enumerate}
\end{athm}

\begin{aproof}
\renewcommand{\u}{\tilde{u}}
Throughout this proof let $c \colon [0,1)\times[0,1) \to [0,\infty)$ satisfy for all $x,y\in[0,1)$ with $x\neq y$ that $c(x,x) = 0$ and
$c(x,y) = (x-y)^{-1}[f(x)-f(y)]$,
let $\u \colon [0,T) \times [-a,a] \allowbreak \to \R$ satisfy for all $t \in [0,T)$, $x \in [-a,a]$ that
$\u(t,-a) = \u(t,a) = 0$, $\u(0,x) = u_0(x)$, 
$ \lim_{s \to T^-} [ \sup_{y\in[-a,a]} \u(s,y) ] = 1$,
and
\begin{equation}\label{cont_2}
\pr[\big]{ \tfrac{\partial}{\partial t} \u}\lrSpace(t,x) = \pr[\big]{ \tfrac{\partial^2}{\partial x^2} \u}\lrSpace(t,x) + f(\u(t,x)) \dc
\end{equation}
let $w \colon [0,T) \times [-a,a] \to \R$ satisfy for all $t \in [0,T)$, $x \in [-a,a]$ that $w(t,x) = u(t,x) - \u(t,x)$,
and
let $v \colon [0,T) \times [-a,a] \to \R$ satisfy for all $t \in[0,T)$, $x\in[-a,a]$ that
$v(t,x) = \pr[]{ \tfrac{\partial }{\partial t} u }(t,x) $.
Observe that \cref{cont_2} implies that for all $t\in[0,T)$, $x\in[-a,a]$ it holds that $w(t,-a) = w(t,a) = w(0,x) = 0$, $\lim_{t\to T^-} w(s,0) = 0$, and
\begin{equation}
\pr[\big]{ \tfrac{\partial}{\partial t} w}\lrSpace(t,x) = \pr[\big]{ \tfrac{\partial^2}{\partial x^2} w}\lrSpace(t,x) + c(u(t,x),\u(t,x)) w(t,x) \dpp
\end{equation}
This ensures that for all $t\in[0,T)$, $x\in[-a,a]$ it holds that $w(t,x) = 0$.
This establishes \cref{lem:quench_props_i0}.
In addition, note that, e.g., Rankin \cite[Theorem 2]{MR1052911} establishes \cref{lem:quench_props_i1}.
Furthermore, observe that the fact that for all $t \in[0,T)$, $x\in[-a,a]$ it holds that
$v(t,x) = \pr[]{ \tfrac{\partial }{\partial t} u }(t,x) $ and \cref{pde2} assure that for all $t\in[0,T)$, $x\in[-a,a]$ it holds that
$v(t,-a) = v(t,a) = 0$, 
\begin{equation}\label{lem:quench_props_eq1}
v(0,x) = \pr[\big]{ \tfrac{\partial^2}{\partial x^2} u_0}\lrSpace(x) + f(u_0(x)) > 0 \dc
\end{equation}
and
\begin{equation}\label{lem:quench_props_eq2}
\pr[\big]{ \tfrac{\partial}{\partial t} v}\lrSpace(t,x) = \pr[\big]{ \tfrac{\partial^2}{\partial x^2} v}\lrSpace(t,x) + f'(u(t,x)) v(t,x) \dpp
\end{equation}
Combining \cref{lem:quench_props_eq1}, \cref{lem:quench_props_eq2}, \cref{pos_init}, \cref{lem:quench_props_i1}, and, e.g., Szarski \cite[Theorem 1]{MR369968} hence assures that for all $t\in[0,T)$, $x\in[-a,a]$ it holds that $v(t,x) \in [0,\infty)$. This and the fact that for all $x\in[-a,a]$ it holds that $u_0(x) \in [0,1)$ establish \cref{lem:quench_props_i2,lem:quench_props_i3}.
\end{aproof}

\begin{definition}[Quenching solution]\label{def:quench}
Assume \cref{setting1} and assume that
\begin{equation}
\textstyle\lim_{t\to T^-} \left[ \sup_{x\in[-a,a]} \pr[\big]{ \tfrac{\partial}{\partial t}u } \lrSpace (t,x) \right] = \infty \dpp
\end{equation}
Then we say that the solution $u$ to \cref{pde2} \emph{quenches}.
Moreover, we refer to $T$ as the \emph{quenching time}.
\end{definition}

\begin{remark}
It is worth noting that \cref{def:quench} excludes the possibility of \emph{quenching in infinite time}.
As such an occurrence is a limiting case of our current studies, the authors feel there is no loss in generality by omitting such a possibility.
\end{remark}

\begin{athm}{lemma}{lem:quench}
Assume \cref{setting1}.
Then it holds that the solution $u$ to \cref{pde2} quenches (cf.\ \cref{def:quench}).
\end{athm}

\begin{aproof}
Throughout this proof let $v \colon [0,T) \times [0,2a] \to \R$ satisfy for all $t\in[0,T)$, $x\in[0,2a]$ that 
\begin{equation}\label{transform1}
v(t,x) = u(t,x-a) - u_0(x-a)
\end{equation}
and let $g_x \in C(\R,\R)$, $x\in[0,2a]$, satisfy for all $x\in[0,2a]$, $w\in [0,1-u_0(x-a))$ that
\begin{equation}\label{transform2}
g_x(w) = f\pr[\big]{ w + u_0(x-a) } + \pr[\big]{ \tfrac{\partial^2}{\partial x^2} u_0}\lrSpace(x-a) \dpp
\end{equation}
Observe that \cref{pde2,transform1,transform2} ensure that for all $t\in[0,T)$, $x\in[0,2a]$ it holds that $v(t,0) = u(t,-a) - u_0(-a) = 0$, $v(t,2a) = u(t,a) - u_0(a) = 0$, $v(0,x) = u(0,x-a) - u_0(x-a) = 0$, and 
\begin{equation}\label{transform3}
\begin{split}
\pr[\big]{ \tfrac{\partial}{\partial t} v}\lrSpace(t,x) 
& = \pr[\big]{ \tfrac{\partial}{\partial t} u}\lrSpace(t,x-a)
- \pr[\big]{ \tfrac{\partial}{\partial t} u_0}\lrSpace(x-a) 
= \pr[\big]{ \tfrac{\partial}{\partial t} u}\lrSpace(t,x-a) \\
& = \pr[\big]{ \tfrac{\partial^2}{\partial x^2} u}\lrSpace(t,x-a) + f(u(t,x-a)) \\
& = \pr[\big]{ \tfrac{\partial^2}{\partial x^2} u}\lrSpace(t,x-a) - \pr[\big]{ \tfrac{\partial^2}{\partial x^2} u_0}\lrSpace(x-a) + f(u(t,x-a)) + \pr[\big]{ \tfrac{\partial^2}{\partial x^2} u_0}\lrSpace(x-a) \\
& = \pr[\big]{ \tfrac{\partial^2}{\partial x^2} v}\lrSpace(t,x) + g_x(v(t,x)) \dpp
\end{split}
\end{equation}
Next, note that \cref{transform2}, \cref{pos_init}, the hypothesis that $f(0) > 0$, the hypothesis that for all $x \in [0,1)$ it holds that $f'(x) > 0$, and the hypothesis that for all $x\in[-a,a]$ it holds that $u_0(x) \in [0,1)$ assure that for all $t\in[0,T)$, $x\in[0,2a]$, $w \in [0,1-u_0(x-a))$ it holds that
\begin{equation}\label{transform4}
g_x(0) = f(u_0(x-a)) + \pr[\big]{ \tfrac{\partial^2}{\partial x^2} u_0}\lrSpace(x-a) > 0
\end{equation}
and
\begin{equation}\label{transform4a}
g_x'(w) = f'\pr[\big]{ w + u_0(x-a) } > 0 \dpp
\end{equation}
In addition, observe that \cref{transform2,f_cond} guarantee that for all $t\in[0,T)$, $x\in[0,2a]$, $w_1,w_2 \in [0,1-u_0(x-a))$, $s\in[0,1]$ it holds that
\begin{equation}\label{transform5}
\abs{g_x(w_1) - g_x(w_2) } = \abs[\big]{ f\pr[\big]{ w_1 + u_0(x-a) } - f\pr[\big]{ w_2 + u_0(x-a) } }
\le L_{w_1,w_2} \abs{w_1-w_2} 
\end{equation}
and
\begin{align}\label{transform6}
g_x\pr[\big]{ sw_1 + (1-s)w_2 } 
& = f\pr[\big]{ sw_1 + (1-s)w_2 + u_0(x-a) } + \pr[\big]{ \tfrac{\partial^2}{\partial x^2} u_0}\lrSpace(x-a) \\
& = f\pr[\big]{ s( w_1+u_0(x-a)) + (1-s)(w_2 + u_0(x-a)) } + \pr[\big]{ \tfrac{\partial^2}{\partial x^2} u_0}\lrSpace(x-a) \nonumber \\
& \le s f\pr[\big]{ w_1 + u_0(x-a) } + (1-s) f\pr[\big]{ w_2 + u_0(x-a) } + \pr[\big]{ \tfrac{\partial^2}{\partial x^2} u_0}\lrSpace(x-a) \nonumber \\
& = sg_x(w_1) + (1-s)g_x(w_2) \dpp \nonumber
\end{align}
Furthermore, note that \cref{transform2}, \cref{f_cond2}, and the assumption that $u_0 \in C^2([-a,a],[0,1))$ yield that for all $t\in[0,T)$, $x\in[0,2a]$ it holds that
\begin{equation}\label{transform7}
\lim_{z\to (1-u_0(x-a))^-} g_x(z) = \lim_{z\to (1-u_0(x-a))^-} \br[\Big]{ f\pr[\big]{ z + u_0(x-a) } + \pr[\big]{ \tfrac{\partial^2}{\partial x^2} u_0}\lrSpace(x-a) }
= \infty
\end{equation}
and
\begin{equation}\label{transform8}
\begin{split}
\int_{0}^{1-u_0(x-a)} g_x(z) \dx z
& = \int_{0}^{1-u_0(x-a)} \br[\Big]{ f\pr[\big]{ z + u_0(x-a) } + \pr[\big]{ \tfrac{\partial^2}{\partial x^2} u_0}\lrSpace(x-a) } \dx z \\
& = \int_0^1 f(z) \dx z + \pr[\big]{ \tfrac{\partial^2}{\partial x^2} u_0}\lrSpace(x-a)
= \infty \dpp
\end{split}
\end{equation}
Combining \cref{transform1,transform2,transform3,%
transform4,transform4a,transform5,transform6,transform7,%
transform8} with Zhou et al.\ \cite[Theorem 3.3]{MR3386569}
(applied for every $x\in[0,2a]$ with $c \with 1-u_0(x-a)$, $a \with 2a$, $u \with v$, $f \with g_x$, $b \with 0$, $T_* \with T$
in the notation of Zhou et al.\ \cite[Theorem 3.3]{MR3386569})
proves that the solution $u$ to \cref{pde2} quenches (cf.\ \cref{def:quench}).
\end{aproof}

\section{Properties of the semidisretized method}\label{sec:3}

In this section, we introduce a semidiscretized approximation of the solution $u$ of \cref{pde2}---such an approach is often referred to as the \emph{method of lines}.
The approximation employed utilizes finite difference operators defined on nonuniform grids.
First, in \cref{sec:3_1} below, we prove results which demonstrate that the semidiscretized solution preserves the observed positivity and monotonicity of the solution $u$ to \cref{pde2}.
Next, in \cref{sec:3_2} below, we introduce a family of discrete weighted $p$-norms and discrete weighted logarithmic norms (see \cref{def:euclid_norm,def:log} below).
We then prove a novel result regarding bounds for the weighted discrete logarithmic norm in the case $p=2$.
Thereafter, in \cref{sec:3_3} below, we provide a rigorous analysis of the convergence of the semidiscrete solution to the solution $u$ of \cref{pde2}.
Moreover, we focus on the effects that the quenching singularity may have on the convergence rate (see \cref{lem:space_conv} below).

\subsection{The semidiscretized method}\label{sec:3_1}

In this subsection, we introduce the employed semidiscretization method in \cref{setting2} below.
This method employs the so-called \emph{central difference} approximation of the one-dimensional Laplace operator on a non-uniform grid.
We prove that the semidiscretized solution is both positive and monotonically increasing on its interval of existence in \cref{lem:space_conv2a} below.
\cref{lem:space_conv2a} in turn, depends on the regularity and representation results presented in \cref{lem:mat_pos_pre} and the matrix positivity result presented in \cref{lem:exp_pos}.
The results in \cref{lem:exp_pos} below are well-known, but we include them for the sake of completeness.

\begin{setting}\label{setting2}
Assume \cref{setting1}, let $\Ts \in (0,\infty)$, $N \in \N$, let $h_0,h_1,\ldots,h_N,x_0,x_1,\ldots,\allowbreak x_{N+1} \allowbreak \in [-a,a]$ satisfy for all $n \in \{0,1,\ldots,N\}$ that $-a = x_0 < x_1 < \ldots < x_{N+1} = a$ and $x_{n+1} - x_n = h_n$, 
let $U \colon [0,\Ts) \to \R^N$ satisfy for all $t \in [0,\Ts)$, $n \in \{1,2,\ldots,N\}$ that $U_n(0) = u_0(x_n)$,
$\lim_{s\to \Ts^-} [ \max_{k\in\{1,2,\ldots,N\}} U_n(s) ] = 1$, and
\begin{equation}\label{eq:semi_dis}
\begin{split}
& \pr[\big]{ \tfrac{{\rm d}}{{\rm d}t} U_n} \lrSpace (t) \\
& \quad = \frac{2 \1_{(1,N]}(n) U_{\max\{1,n-1\}}(t)}{h_{n-1}(h_{n-1} + h_n)} - \frac{2 U_n(t) }{h_{n-1}h_n} + \frac{2 \1_{[1,N)}(n) U_{\min\{N,n+1\}}(t)}{h_n(h_{n-1}+h_n)} + f(U_n(t)) \dc
\end{split}
\end{equation}
let $F \colon \R^N \to \R^N$ satisfy for all $X \in \R^N$, $n \in \{1,2,\ldots,N\}$ that $(F(X))_n = f(X_n)$, 
and let $A = (A_{i,j})_{i,j \in\{1,2,\ldots,N\}} \in \R^{N\times N}$ satisfy for all $i,j,n \in \{1,2,\ldots,N\}$, $k \in \{1,2,\ldots,N-1\}$ with $\abs{i-j} \in \{2,3,\ldots,N-1\}$ that $A_{i,j} = 0$,
\begin{equation}\label{eq:mat_A}
A_{k+1,k} = \frac{2}{h_k(h_k + h_{k+1})} \dc
\quad
A_{n,n} = \frac{-2}{h_{n-1}h_n} \dc
\quad \text{and} \quad
A_{k,k+1} = \frac{2}{h_k(h_{k-1}+h_k)} \dpp
\end{equation}
\end{setting}

\begin{remark}
Note that, just as in \cref{setting1}, we have assumed that $\Ts \in (0,\infty)$ is the quenching time for our semidiscrete approximation to the solution to \cref{pde2} (cf.\ \cref{def:quench}).
Observe that this assumption has implicitly added additional conditions on the value $a\in(0,\infty)$ (as $a$ must be chosen so that the true solution \emph{and} the semidiscrete solution will exhibit quenching behavior).
For an example of how one may study the relationship between $T$ and $\Ts$, interested readers may refer to, e.g., Nabongo and Boni \cite{MR2490440}.
\end{remark}

\begin{definition}[Matrix exponential]\label{def:mat_exp}
We denote by $\exp \colon \bigcup_{d\in\N} \! \R^d \to \bigcup_{d\in\N} \! \R^d$ the function which satisfies for all $d \in \N$, $B\in\R^{d\times d}$ that
$\exp(B) = \sum_{k=0}^{\infty} (\nicefrac{1}{k!}) B^k$.
\end{definition}

\begin{definition}[Matrix comparison operations]\label{def:wedge}
Let $m,n\in\N$ and let $B,C \in \R^{m\times n}$ satisfy for all $i \in \{1,2,\ldots,m\}$, $j\in\{1,2,\ldots,n\}$ that $B_{i,j} \ge C_{i,j}$ (similarly, $B_{i,j} > C_{i,j}$).
Then we say that $B \ge C$ (similarly $B > C$).
Moreover, if for all $i \in \{1,2,\ldots,m\}$, $j\in\{1,2,\ldots,n\}$ it holds that $C_{i,j} = 0$, then we say that $B$ is \emph{nonnegative} (similarly, $B$ is \emph{positive}).
\end{definition}

\begin{athm}{proposition}{lem:mat_compare}
Let $d_1,d_2,d_3 \in \N$, let $B_1 \in \R^{d_1\times d_2}$ be nonsingular and nonnegative, and let $B_2 \in \R^{d_2\times d_3}$ be nonsingular and positive (cf.\ \cref{def:wedge}).
Then it holds that $B_1B_2$ is positive.
\end{athm}

\begin{aproof}
Note that the result is an immediate consequence of \cref{def:wedge}.
\end{aproof}

\begin{definition}[Identity matrix]\label{def:identityMatrix}
Let $d\in\N$. Then we denote by $\idMatrix_{d}\in \R^{d\times d}$ the identity matrix in $\R^{d\times d}$.
\end{definition}

\begin{athm}{lemma}{lem:exp_pos}
Assume \cref{setting2}.
Then it holds for all $t\in[0,\infty)$ that $\exp(tA)$ is positive (cf.\ \cref{def:mat_exp,def:wedge}).
\end{athm}

\begin{aproof}
Throughout this proof let $\alpha = \nicefrac{2}{(\min_{n\in\{1,2,\ldots,N\}} h_{n-1}h_n)}$ and let $B \in \R^{N\times N}$ satisfy $B = A + \alpha \idMatrix_N$ (cf.\ \cref{def:identityMatrix}).
Note that the fact that for all $C \in \R^{N\times N}$ it holds that $CI_N = I_NC$ guarantees that for all $t\in[0,\infty)$ it holds that
\begin{equation}\label{lem:exp_pos_1}
\exp(tA) = \exp\pr[\big]{ t(B - \alpha \idMatrix_N) } = \exp(-\alpha t\idMatrix_N) \exp(tB)
\end{equation}
(cf.\ \cref{def:mat_exp}).
Next, observe that \cref{eq:mat_A} ensures that $B$ is nonnegative (cf.\ \cref{def:wedge}).
This assures that for all $t\in[0,\infty)$ it holds that $\exp(tB)$ is nonnegative.
Combining this, \cref{lem:exp_pos_1}, the fact that for all $t\in[0,\infty)$ it holds that $\exp(-\alpha t\idMatrix_N)$ is positive, and \cref{lem:mat_compare} proves that for all $t\in[0,\infty)$ it holds that $\exp(tA)$ is positive.
\end{aproof}

\begin{athm}{lemma}{lem:mat_pos_pre}
Assume \cref{setting2}.
Then 
\begin{enumerate}[label=(\roman*)]
\item\label{lem:space_conv2a_i0} it holds that $U \in C^1([0,\Ts),\R^N)$,
\item\label{lem:space_conv2a_i1} it holds for all $t \in [0,\Ts)$ that
\begin{equation}\label{eq:semi_dis2}
\pr[\big]{ \tfrac{{\rm d}}{{\rm d}t} U} \lrSpace (t)
= AU(t) + F(U(t)) \dc
\end{equation}
\item\label{lem:space_conv2a_i2} it holds for all $t \in [0,\Ts)$ that
\begin{equation}\label{eq:semi_dis2a}
U(t) = \exp(tA) U(0) + \int_0^t \exp\pr[\big]{ (t-s) A } F(U(s)) \dx s
\end{equation}
\end{enumerate}
(cf.\ \cref{def:mat_exp,def:wedge}).
\end{athm}

\begin{aproof}
First, note that, e.g., Rankin \cite[Theorem 2]{MR1052911} establishes \cref{lem:space_conv2a_i0}.
Next,
observe that \cref{eq:semi_dis2}
is an immediate consequence of combining \cref{lem:space_conv2a_i0} with \cref{eq:semi_dis,eq:mat_A}.
This establishes \cref{lem:space_conv2a_i1}.
In addition, note that solving \cref{eq:semi_dis2} via standard integrating factor techniques yields \cref{eq:semi_dis2a}.
This establishes \cref{lem:space_conv2a_i2}.
\end{aproof}

\begin{athm}{lemma}{lem:space_conv2a}
Assume \cref{setting2} and assume that $AU(0) + F(U(0))$ is positive (cf.\ \cref{def:wedge}).
Then 
\begin{enumerate}[label=(\roman*)]
\item\label{lem:space_conv2a_i3a} it holds for all $t \in [0,\Ts)$ that $U(t)$ is positive,
\item\label{lem:space_conv2a_i3} it holds for all $t_1, t_2 \in [0,\Ts)$ with $t_1 \le t_2$ that $U(t_1) \le U(t_2)$, and
\item\label{lem:mat_pos_pre_i4} it holds for all $t\in[0,\Ts)$ that $AU(t) + F(U(t))$ is positive
\end{enumerate}
(cf.\ \cref{def:mat_exp}).
\end{athm}

\begin{aproof}
Throughout this proof let $V \colon [0,\Ts) \to \R^N$ satisfy for all $t\in[0,\Ts)$ that 
\begin{equation}\label{eq:der_U_fn}
V(t) = \pr[\big]{ \tfrac{{\rm d}}{{\rm d}t} U} \lrSpace (t) \dpp
\end{equation}
Observe that \cref{eq:der_U_fn}, \cref{lem:space_conv2a_i0,lem:space_conv2a_i1} in \cref{lem:mat_pos_pre}, and the assumption that $AU(0) + F(U(0))$ is positive demonstrate that for all $t\in[0,\Ts)$ it holds that
\begin{equation}\label{eq:der_U_fn1}
V(0) = AU(0) + F(U(0)) > 0
\end{equation}
and
\begin{equation}\label{eq:der_U_fn2}
\pr[\big]{ \tfrac{{\rm d}}{{\rm d}t} V} \lrSpace (t) = A V(t) + F'(U(t)) V(t) \dpp
\end{equation}
Combining \cref{eq:der_U_fn1}, \cref{eq:der_U_fn2}, the assumption that $AU(0) + F(U(0))$ is positive, \cref{lem:space_conv2a_i0} in \cref{lem:mat_pos_pre}, and, e.g., Szarski \cite[Theorem 1]{MR369968} hence assures that for all $n\in\{1,2,\ldots,N\}$, $t\in[0,\Ts)$ it holds that $V_n(t) \in [0,\infty)$. This and the fact that for all $n\in\{1,2,\ldots,N\}$ it holds that $U_n(0) \in [0,1)$ establish \cref{lem:space_conv2a_i3a,lem:space_conv2a_i3}.
We now prove \cref{lem:mat_pos_pre_i4} by transfinite induction on $t\in[0,\Ts)$ (cf., e.g., Chao \cite{MR1560239}).
Note that the base case $t\in[0,0]$ holds due to the assumption that $AU(0) + F(U(0))$ is positive.
This establishes that $AU(t) + F(U(t))$ is positive in the base case $t\in[0,0]$.
For the transfinite induction step $[0,\Ts) \supseteq [0,T_*) \ni t \to t \in [0,T_*] \subseteq [0,\Ts)$, let $T_*\in[0,\Ts)$ and assume that for every $t \in [0,T_*)$ it holds that 
\begin{equation}\label{eq:tran_ind}
AU(t) + F(U(t)) > 0 \dpp
\end{equation}
Observe that \cref{lem:space_conv2a_i2} in \cref{lem:mat_pos_pre} and integration by parts guarantee that 
\begin{align}\label{eq:tran_ind1}
& AU(T_*) + F(U(T_*)) 
= A \left[ \exp(T_*A) U(0) + \int_0^{T_*} \exp\pr[\big]{ (T_*-s) A } F(U(s)) \dx s \right] + F(U(T_*)) \nonumber \\
& \quad = \exp(T_*A) \left[ AU(0) + \int_0^{T_*} A \exp(-sA) F(U(s)) \dx s + \exp(-T_*A) F(U(T_*)) \right] \nonumber \\
& \quad = \exp(T_*A) \Biggl[ AU(0) + \Bigl( F(U(0)) - \exp(-T_*A) F(U(T_*)) \Bigr) \\
& \quad \qquad + \int_0^{T_*} \exp(-sA) F'(U(s)) \pr[\big]{ \tfrac{{\rm d}}{{\rm d}s} U} \lrSpace (s) \dx s + \exp(-T_*A) F(U(T_*)) \Biggr] \nonumber \\
& \quad = \exp(T_*A) \br[\big]{ AU(0) + F(U(0)) } + \int_0^{T_*} \exp\pr[\big]{(T_*-s)A} F'(U(s)) \pr[\big]{ \tfrac{{\rm d}}{{\rm d}s} U} \lrSpace (s) \dx s \nonumber
\end{align}
(cf.\ \cref{def:mat_exp}).
Moreover, note that \cref{lem:mat_compare}, \cref{lem:exp_pos}, \cref{lem:space_conv2a_i3a}, \cref{eq:tran_ind}, and the assumption that for all $x\in[0,1)$ it holds that $f'(x) > 0$ proves that
\begin{equation}\label{eq:tran_ind2}
\begin{split}
& \int_0^{T_*} \exp\pr[\big]{(t-s)A} F'(U(s)) \pr[\big]{ \tfrac{{\rm d}}{{\rm d}s} U} \lrSpace (s) \dx s \\
& \qquad = \int_0^{T_*} \exp\pr[\big]{(t-s)A} F'(U(s)) \br[\big]{ AU(s) + F(U(s)) } \dx s
> 0 \dpp
\end{split}
\end{equation}
Combining \cref{eq:tran_ind1}, \cref{eq:tran_ind2}, \cref{lem:mat_compare}, \cref{lem:exp_pos}, and the assumption that $AU(0) + F(U(0))$ is positive hence shows that  
$AU(T_*) + F(U(T_*))$ is positive.
Transfinite induction hence establishes \cref{lem:mat_pos_pre_i4}.
\end{aproof}

\subsection{Weighted norms and logarithmic norms}\label{sec:3_2}

In this subsection, we introduce the discrete weighted $p$-norms and their associated discrete weighted logarithmic norms.
The weighted $p$-norms in \cref{def:euclid_norm} below may be understood as approximations of the standard $L^p$-norms via nonuniform quadrature methods.
We present these results for any general sequence of weights (denoted by $\fh$ below).
However, as such studies are not of immediate interest, we do not present any details regarding the convergence of the norms in \cref{def:euclid_norm} below to their continuous counterparts.
Moreover, we do not prove here that \cref{def:euclid_norm} is indeed a norm, but this fact should be obvious.

The primary obstacle when employing nonuniform meshes in the semidiscretization procedure is the fact that the resulting coefficient matrix may no longer be symmetric (cf., e.g., \cref{eq:mat_A}).
Thus, any ensuing norm estimates become more difficult to provide sharp bounds for.
In order to circumvent this issue, we employ the logarithmic norm; however, despite its namesake, the logarithmic norm is not in fact a norm.
The logarithmic norm has a long history of being studied and employed in numerical analysis; interested readers may refer to, e.g., \cite{MR408227,MR2265579,Sheng2} for more details.
Herein, we are primarily interested in the use of the associated discrete weighted logarithmic norms for estimating bounds related to the matrix $A$ in \cref{eq:mat_A}.
As such, the main result of this subsection is \cref{lem:A_log} below.
The proof of \cref{lem:A_log} is in turn based on the result in \cref{lem:log_rep}. 

\begin{definition}[Weighted $p$-norm]\label{def:euclid_norm}
For every $d \in \N$, $\fh_0,\fh_1,\ldots,\fh_d \in (0,\infty)$, $\fh = (\fh_0,\fh_1,\ldots,\fh_d) \in \R^{d+1}$
we denote by $\norm{\cdot}_{\fh,p} \colon \R^d \to [0,\infty)$, $p\in[1,\infty]$, the functions which satisfy for all $p\in[1,\infty)$, 
$x = (x_1,x_2,\dots,x_d) \in \R^d$ that $\norm{x}_{\fh,p} = \br{ \sum_{k=1}^d \tfrac{1}{2}(\fh_{k-1}+\fh_k) \abs{x_k}^p }^{\nicefrac{1}{p}}$ and $\norm{x}_{\fh,p} = \max_{k\in\{1,2,\ldots,d\}} \abs{x_k}$.
\end{definition}

\begin{definition}[Weighted $p$-logarithmic norm]\label{def:log}
For every $d \in \N$, $\fh_0,\fh_1,\ldots,\fh_d \in (0,\infty)$, $\fh = (\fh_0,\fh_1,\ldots,\fh_d)\in\R^{d+1}$
we denote by $\logNorm_{\fh,p} \colon \R^{d\times d} \to \R$, $p\in[1,\infty]$, the functions which satisfy for all $p\in[1,\infty]$, $B \in \R^{d\times d}$ that
\begin{equation}
\logNorm_{\fh,p}(B) = \sup_{ \substack{ v \in \R^d \\ \norm{v}_{\fh,p} \neq 0 } } \left[ \lim_{t\to 0^+} \frac{ \norm[]{ (\idMatrix_d + t B ) v }_{\fh,p} - \norm{v}_{\fh,p} }{ t \norm{v}_{\fh,p} } \right]
\end{equation}
(cf.\ \cref{def:euclid_norm,def:identityMatrix}).
\end{definition}

\begin{definition}[Matrix transpose]\label{def:Transpose}
Let $m, n \in \N$, $B \in \R^{m \times n}$. Then we denote by $B^* \in \R^{n \times m}$ the transpose of $B$.
\end{definition}

\begin{athm}{lemma}{lem:log_rep}
Let $d\in\N$, $\fh_0,\fh_1,\ldots,\fh_d \in (0,\infty)$, $\fh = (\fh_0,\fh_1,\ldots,\fh_d)\in\R^{d+1}$ and let $H \in \R^{d\times d}$ be the matrix which satisfies for all $i,j\in\{1,2,\ldots,d\}$ with $i\neq j$ that
$H_{i,i} = \sqrt{\nicefrac{\fh_{i-1}+\fh_i}{2}}$ and $H_{i,j} = 0$.
Then it holds for all $B \in\R^{d\times d}$ that
\begin{align}
& \logNorm_{\fh,2}(B) \\
& = \max\cu[\big]{ \lambda \in \R \colon \pr[\big]{ \exists \, v \in \R^d \text{ with } \norm{v}_{\fh,2} \neq 0 \text{ and } [(HBH^{-1})^* + (HBH^{-1})]v = 2 \lambda v } } \dpp \nonumber
\end{align}
(cf.\ \cref{def:log,def:Transpose,def:euclid_norm}).
\end{athm}

\begin{aproof}
Throughout this proof 
let $\vt{\cdot,\cdot} \colon \R^d\times \R^d \to \R$ be the function which satisfies for all $v,w\in\R^d$ that $\vt{v,w} = \sum_{k=1}^d v_kw_k$.
Note that for all $v \in \R^d$ it holds that
\begin{equation}\label{eq:inner_to_norm}
\vt{Hv,Hv} = \SmallSum_{k=1}^d (Hv)_k(Hv)_k 
= \SmallSum_{k=1}^d \tfrac{1}{2}(\fh_{k-1}+\fh_k) \abs{v_k}^2 = \norm{v}_{\fh,2}^2
\end{equation}
(cf.\ \cref{def:euclid_norm}).
This implies that for all $t\in\R$, $v\in\R^d$, $B\in\R^{d\times d}$ it holds that
\begin{equation}\label{log_est1}
\begin{split}
\frac{ \norm[]{ (\idMatrix_d + t B ) v }_{\fh,p} - \norm{v}_{\fh,p} }{ t \norm{v}_{\fh,p} }
& = \frac{ \norm[]{ (\idMatrix_d + t B ) v }_{\fh,p}^2 - \norm{v}_{\fh,p}^2 }{ t \norm{v}_{\fh,p} [ \norm[]{ (\idMatrix_d + t B ) v }_{\fh,p} + \norm{v}_{\fh,p} ] } \\
& = 
\frac{ \vt{ H (\idMatrix_d + t B ) v , H (\idMatrix_d + t B ) v } - \vt{Hv,Hv} }{ t \norm{v}_{\fh,p} \br[\big]{ \norm[]{ (\idMatrix_d + t B ) v }_{\fh,p} + \norm{v}_{\fh,p} } }
\end{split}
\end{equation}
(cf.\ \cref{def:identityMatrix}).
Next, observe that for all $t\in\R$, $v\in\R^d$, $B\in\R^{d\times d}$ it holds that
\begin{equation}
\vt{ H (\idMatrix_d + t B ) v , H (\idMatrix_d + t B ) v } - \vt{Hv,Hv}
= t \vt{ HBv,Hv } + t \vt{ Hv, HBv } \dpp
\end{equation}
Combining this and \cref{log_est1} ensures that for all $v\in\R^d$, $B\in\R^{d\times d}$ it holds that
\begin{equation}\label{log_est2}
\begin{split}
\lim_{t\to 0} \frac{ \norm[]{ (\idMatrix_d + t B ) v }_{\fh,p} - \norm{v}_{\fh,p} }{ t \norm{v}_{\fh,p} }
& = \lim_{t\to 0} \frac{ \vt{ H (\idMatrix_d + t B ) v , H (\idMatrix_d + t B ) v } - \vt{Hv,Hv} }{ t \norm{v}_{\fh,p} \br[\big]{ \norm[]{ (\idMatrix_d + t B ) v }_{\fh,p} + \norm{v}_{\fh,p} } } \\
& = \lim_{t\to 0} \frac{ t \vt{ HBv,Hv } + t \vt{ Hv, HBv } }{ t \norm{v}_{\fh,p} \br[\big]{ \norm[]{ (\idMatrix_d + t B ) v }_{\fh,p} + \norm{v}_{\fh,p} } } \\
& = \frac{ \vt{ HBv,Hv } + \vt{ Hv, HBv } }{ 2\norm{v}_{\fh,p}^2 } \dpp
\end{split}
\end{equation}
Moreover, note that the assumption that $\fh_0,\fh_1,\ldots,\fh_d \in (0,\infty)$ guarantees that for all $v \in \R^d$, $B\in\R^{d\times d}$ it holds that
\begin{equation}
\begin{split}
\vt{ HBv,Hv } + \vt{ Hv, HBv } 
& = \vt[\big]{ (HBH^{-1}) Hv , Hv } + \vt[\big]{ Hv , (HBH^{-1}) Hv } \\
& = \vt[\big]{ Hv , (HBH^{-1})^* Hv } + \vt[\big]{ Hv , (HBH^{-1}) Hv } \\
& = \vt[\big]{ Hv , [(HBH^{-1})^* + (HBH^{-1})] Hv }
\end{split}
\end{equation}
(cf.\ \cref{def:Transpose}).
This and \cref{log_est2} assure that for all $B\in\R^{d\times d}$ it holds that
\begin{equation}\label{log_est3}
\begin{split}
\logNorm_{\fh,2}(B)
& = \sup_{ \substack{ v \in \R^d \\ \norm{v}_{\fh,p} \neq 0 } } \left[ \frac{ \vt{ HBv,Hv } + \vt{ Hv, HBv } }{ 2\norm{v}_{\fh,p}^2 } \right] \\
& =
\sup_{ \substack{ v \in \R^d \\ \norm{v}_{\fh,p} \neq 0 } } \left[ \frac{ \vt[\big]{ Hv , \tfrac{1}{2}[(HBH^{-1})^* + (HBH^{-1})] Hv } }{ \vt{Hv,Hv} } \right] \\
& =
\sup_{ \substack{ w \in \R^d \\ \norm{w}_{\fh,p} \neq 0 } } \left[ \frac{ \vt[\big]{ w , \tfrac{1}{2}[(HBH^{-1})^* + (HBH^{-1})] w } }{ \vt{w,w} } \right] 
\end{split}
\end{equation}
(cf.\ \cref{def:log}).
Combining \cref{log_est3} and and the Rayleigh quotient theorem (cf., e.g., Driver \cite[Theorem A.26]{driverfunctional})
hence proves that for all $B\in\R^{d\times d}$ it holds that
\begin{align}
& \logNorm_{\fh,2}(B) \\
& = \max\cu[\big]{ \lambda \in \R \colon \pr[\big]{ \exists \, v \in \R^d \text{ with } \norm{v}_{\fh,2} \neq 0 \text{ and } [(HBH^{-1})^* + (HBH^{-1})]v = 2 \lambda v } } \dpp \nonumber
\end{align}
\end{aproof}

\begin{athm}{lemma}{lem:A_log}
Assume \cref{setting2} and let $\fh \in \R^{N+1}$ satisfy $\fh = (h_0,h_1,\ldots,h_N)$.
Then it holds that $\logNorm_{\fh,2}(A) \in (-\infty,0)$ (cf.\ \cref{def:log}).
\end{athm}

\begin{aproof}
Throughout this proof let $H \in \R^{N\times N}$ be the matrix which satisfies for all $i,j\in\{1,2,\ldots,N\}$ with $i\neq j$ that
$H_{i,i} = \sqrt{\nicefrac{h_{k-1}+h_k}{2}}$ and $H_{i,j} = 0$.
Note that for all $B \in \R^{N\times N}$ it holds that
\begin{equation}
(HBH^{-1})^* + (HBH^{-1}) = H^{-1}B^*H + HBH^{-1} = H \br[\big]{ B + H^{-2} B^* H^2 } H^{-1} 
\end{equation}
(cf.\ \cref{def:Transpose}).
Combining this and \cref{lem:log_rep} assures that
\begin{align}\label{disk1}
& \logNorm_{\fh,2}(A) \\
& = \max\cu[\big]{ \lambda \in \R \colon \pr[\big]{ \exists \, v \in \R^d \text{ with } \norm{v}_{\fh,2} \neq 0 \text{ and } [(HAH^{-1})^* + (HAH^{-1})]v = 2 \lambda v } } \nonumber \\
& = \max\cu[\big]{ \lambda \in \R \colon \pr[\big]{ \exists \, v \in \R^d \text{ with } \norm{v}_{\fh,2} \neq 0 \text{ and } [ A + H^{-2} A^* H^2 ]v = 2 \lambda v } } \nonumber
\end{align}
(cf.\ \cref{def:log}).
Next, observe that for all $i,j,n \in \{1,2,\ldots,N\}$, $k\in\{1,2,\ldots,N-1\}$ with $\abs{i-j} \in \{2,3,\ldots,N-1\}$ it holds that $(H^{-2} A^* H)_{i,j} = 0$, $(H^{-2} A^* H)_{n,n} = A_{n,n}$,
\begin{equation}\label{disk1a}
(H^{-2} A^* H)_{k+1,k} = \tfrac{h_{k-1}+h_{k}}{h_{k}+h_{k+1}} A_{k,k+1} \dc
\qquad \text{and} \qquad
(H^{-2} A^* H)_{k,k+1} = \tfrac{h_{k}+h_{k+1}}{h_{k-1}+h_{k}} A_{k+1,k} \dpp
\end{equation}
In addition, note that the assumption that $h_0,h_1,\ldots,h_N \in (0,\infty)$ demonstrates that for all $n \in \{1,2,\ldots,N\}$ it holds that
\begin{align}\label{disk2}
& \1_{(1,N]}(n) A_{\max\{n-1,1\},n} + A_{n,n} + \1_{[1,N)}(n) A_{n,\min\{n-1,N\}} \nonumber \\
& \quad = \1_{(1,N]}(n) \frac{2}{h_{n-1}(h_{n-1} + h_{n})} + \frac{-2}{h_{n-1}h_n} + \1_{[1,N)}(n) \frac{2}{h_n(h_{n-1}+h_n)} \nonumber \\
& \quad = 2 \left[ \frac{ \1_{(1,N]}(n) h_n - (h_{n-1}+h_n) + \1_{[1,N)}(n)h_{n-1} }{ h_{n-1}h_n(h_{n-1}+h_n)} \right] \\
& \quad = 2 \left[ \frac{ (\1_{(1,N]}(n)-1) h_n + (\1_{[1,N)}(n)-1)h_{n-1} }{ h_{n-1}h_n(h_{n-1}+h_n)} \right]
= -2 \left[ \frac{ \1_{\{1\}}(n) h_n + \1_{\{N\}}(n) h_{n-1} }{ h_{n-1}h_n(h_{n-1}+h_n)} \right] \dpp \nonumber
\end{align}
Moreover, observe that \cref{disk1a} the fact that $h_0,h_1,\ldots,h_N \in (0,\infty)$ show that for all $n \in \{1,2,\ldots,N\}$ it holds that
\begin{align}\label{disk3}
& \1_{(1,N]}(n) (H^{-2} A^* H)_{\max\{n-1,1\},n} + (H^{-2} A^* H)_{n,n} + \1_{[1,N)}(n) (H^{-2} A^* H)_{n,\min\{n-1,N\}} \nonumber \\
& \quad = \1_{(1,N]}(n) \frac{2}{h_{n-1}(h_{n-1}+h_{n})} + \frac{-2}{h_{n-1}h_n} + \1_{[1,N)}(n) \frac{2}{h_n(h_{n-1}+h_n)} \nonumber \\
& \quad = 2 \left[ \frac{ \1_{(1,N]}(n) h_n - (h_{n-1}+h_n) + \1_{[1,N)}(n)h_{n-1} }{ h_{n-1}h_n(h_{n-1}+h_n)} \right] \\
& \quad = 2 \left[ \frac{ (\1_{(1,N]}(n)-1) h_n + (\1_{[1,N)}(n)-1)h_{n-1} }{ h_{n-1}h_n(h_{n-1}+h_n)} \right]
= -2 \left[ \frac{ \1_{\{1\}}(n) h_n + \1_{\{N\}}(n) h_{n-1} }{ h_{n-1}h_n(h_{n-1}+h_n)} \right] \dpp \nonumber
\end{align}
Combining \cref{disk1,disk2,disk3} with the Ger{\v s}grin criterion (cf., e.g., Iserles \cite[Lemma 8.3]{Iserles})
hence proves that
$\logNorm_{\fh,2}(A) \in (-\infty,0)$.
\end{aproof}

\subsection{Convergence analysis of the semidiscretized method}\label{sec:3_3}

In this subsection we analyze the convergence of the semidiscrete approximation in \cref{eq:semi_dis} to the solution $u$ to \cref{pde2} via the discrete weighted $2$-norm introduced in the previous subsection (cf.\ \cref{def:euclid_norm}).
The aforementioned convergence result is contained in \cref{lem:space_conv}.
Note that the result in \cref{lem:space_conv} emphasizes the dependence of the convergence rate on the regularity of the solution $u$ to \cref{pde2} and the given Lipschitz function in \cref{setting1} (cf.\ \cref{f_cond}).
The proof of \cref{lem:space_conv} utilizes the preliminary result given in \cref{lem:space_conv_pre} below.
\cref{lem:space_conv_pre} provides a representation of the exact error incurred by approximating the one-dimensional Laplace operator with the nonuniform finite difference method introduced in \cref{setting2}.

Note that throughout the remainder of this article we assume without loss of generality that $\Ts \le T$.

\begin{athm}{lemma}{lem:space_conv_pre}
Assume \cref{setting2}.
Then it holds for all $t\in[0,\Ts)$, $n \in \{1,2,\ldots,N\}$ that
\begin{align}\label{eq:space_taylor0}
& \1_{(1,N]}(n) A_{n,n-1} u(t,x_{n-1})+ A_{n,n} u(t,x_n) + \1_{[1,N)}(n) A_{n,n+1} u(t,x_{n+1} ) - \pr[\big]{ \tfrac{\partial^2}{\partial x^2} u}\lrSpace(t,x_n) \nonumber \\
& \quad = - \int_0^{h_{n-1}} \frac{\pr[\big]{ \tfrac{\partial^3}{\partial x^3} u}\lrSpace(t,s)}{h_{n-1}(h_{n-1}+h_n)} (h_{n-1}-s)^2 \dx s +  \int_0^{h_{n}} \frac{\pr[\big]{ \tfrac{\partial^3}{\partial x^3} u}\lrSpace(t,s)}{h_{n}(h_{n-1}+h_n)} (h_{n}-s)^2 \dx s \dpp
\end{align}
\end{athm}

\begin{aproof}
Note that Taylor's theorem (cf., e.g., Cartan et al.\ \cite[Theorem 5.6.3]{cartan2017differential}) ensures that for all $t \in [0,\Ts)$, $n\in\{1,2,\ldots,N\}$ it holds that
\begin{align}\label{eq:space_taylor1}
& u(t,x_{n-1}) = u(t,x_n - h_{n-1}) \\
& = u(t,x_n) - h_{n-1} \pr[\big]{ \tfrac{\partial}{\partial x} u}\lrSpace(t,x_n) + \tfrac{(h_{n-1})^2}{2} \pr[\big]{ \tfrac{\partial^2}{\partial x^2} u}\lrSpace(t,x_n) - \int_0^{h_{n-1}} \frac{\pr[\big]{ \tfrac{\partial^3}{\partial x^3} u}\lrSpace(t,s)}{2} (h_{n-1}-s)^2 \dx s \nonumber
\end{align}
and
\begin{align}\label{eq:space_taylor2}
& u(t,x_{n+1}) = u(t,x_n + h_{n}) \\
& = u(t,x_n) + h_{n} \pr[\big]{ \tfrac{\partial}{\partial x} u}\lrSpace(t,x_n) + \tfrac{(h_{n})^2}{2} \pr[\big]{ \tfrac{\partial^2}{\partial x^2} u}\lrSpace(t,x_n) + \int_0^{h_{n}} \frac{\pr[\big]{ \tfrac{\partial^3}{\partial x^3} u}\lrSpace(t,s)}{2} (h_{n}-s)^2 \dx s \dpp \nonumber
\end{align}
Combining \cref{eq:space_taylor1}, \cref{eq:space_taylor2}, and the assumption that for all $t\in[0,T]$ it holds that $u(t,-a)=u(t,a) = 0$ hence establishes \cref{eq:space_taylor0}.
\end{aproof}

\begin{athm}{lemma}{lem:space_conv}
Assume \cref{setting2}, let $\U \colon [0,\Ts) \allowbreak \to \R^N$ satisfy for all $t\in[0,\Ts)$, $n\in\{1,2,\ldots,N\}$ that $\U_n(t) = u(t,x_n)$, let $\fh \in \R^{N+1}$ satisfy $\fh = (h_0,h_1,\ldots,h_N)$, and let $\LL \colon [0,\Ts) \to [0,\infty)$ satisfy for all $t\in[0,\Ts)$ that $\LL_t = \max_{n\in\{1,2,\ldots,N\}} L_{u(t,x_n),U_n(t)}$.
Then for all $t \in [0,\Ts)$
it holds that
\begin{equation}\label{eq:sd_conv_bd_eq}
\norm{ \U(t) - U(t) }_{\fh,2} \le \sqrt{2a} \exp(\Ts \LL_t) \left[ \int_0^t \max_{n\in\{1,2,\ldots,N\}} \abs*{ \int_0^{h_n} \pr[\big]{ \tfrac{\partial^3}{\partial x^3} u}\lrSpace(w,s) \dx s } \dx w \right]
\end{equation}
(cf.\ \cref{def:euclid_norm}).
\end{athm}

\begin{aproof}
Throughout this proof for every $v \colon [0,\Ts) \to \R$, $t\in[0,\Ts)$ let $D_tv(t) \in [-\infty,\infty]$ satisfy
$D_t v(t) = \limsup_{\varepsilon \to 0^+} \varepsilon^{-1} [ v(t+\varepsilon) - v(t) ] $, let $\V \colon [0,\Ts) \to \R^N$ satisfy for all $t\in[0,\Ts)$ that $\V(t) = \U(t) - U(t)$,
and let $\CC \colon [0,\Ts) \to \R^N$ satisfy for all $n \in \{1,2,\ldots,N\}$ that
\begin{equation}
\CC_n(t) = - \int_0^{h_{n-1}} \frac{\pr[\big]{ \tfrac{\partial^3}{\partial x^3} u}\lrSpace(t,s)}{h_{n-1}(h_{n-1}+h_n)} (h_{n-1}-s)^2 \dx s +  \int_0^{h_{n}} \frac{\pr[\big]{ \tfrac{\partial^3}{\partial x^3} u}\lrSpace(t,s)}{h_{n}(h_{n-1}+h_n)} (h_{n}-s)^2 \dx s \dpp
\end{equation}
Note that \cref{pde2} and \cref{lem:space_conv_pre} ensure that for all $t\in[0,\Ts)$ it holds that
\begin{equation}
\pr[\big]{ \tfrac{{\rm d}}{{\rm d}t} \U}\lrSpace(t) = A\U(t) + F(\U(t)) + \CC(t) \dpp
\end{equation}
This, Taylor's theorem (cf., e.g., Cartan et al.\ \cite[Theorem 5.6.3]{cartan2017differential}), the triangle inequality, and \cref{f_cond} imply that for all $t\in[0,\Ts)$ it holds that
\begin{align}
& D_t \norm{ \V(t) }_{\fh,2}
= \limsup_{\varepsilon\to 0^+} \varepsilon^{-1} \br[\big]{ \norm{ \V(t+\varepsilon) }_{\fh,2} - \norm{ \V(t) }_{\fh,2} } \nonumber \\
& \quad = \limsup_{\varepsilon\to 0^+} \varepsilon^{-1} \br[\Big]{ \norm[\big]{ \V(t) + \varepsilon \pr[\big]{ \tfrac{{\rm d}}{{\rm d}t} \V}\lrSpace(t) }_{\fh,2} - \norm{ \V(t) }_{\fh,2} } \\
& \quad = \limsup_{\varepsilon\to 0^+} \varepsilon^{-1} \br[\Big]{ \norm[\big]{ \V(t) + \varepsilon \pr[\big]{ A\V(t) + F(\U(t)) - F(U(t)) + \CC(t) } }_{\fh,2} - \norm{ \V(t) }_{\fh,2} } \nonumber \\
& \quad \le \norm{\CC(t)}_{\fh,2} + \norm{ F(\U(t)) - F(U(t)) }_{\fh,2} + \limsup_{\varepsilon\to 0^+} \varepsilon^{-1} \br[\Big]{ \norm[\big]{ \V(t) + \varepsilon A\V(t) }_{\fh,2} - \norm{ \V(t) }_{\fh,2} } \nonumber \\
& \quad \le \norm{\CC(t)}_{\fh,2} + \LL_t \norm{ \V(t) }_{\fh,2} + \logNorm_{\fh,2}(A) \norm{ \V(t) }_{\fh,2}
= \norm{\CC(t)}_{\fh,2} + \pr[\big]{ \LL_t + \logNorm_{\fh,2}(A) } \norm{ \V(t) }_{\fh,2}\nonumber 
\end{align}
(cf.\ \cref{def:log,def:euclid_norm}).
This, the fact that for all $n\in\{1,2,\ldots,N\}$ it holds that $u(0,x_n) = U_n(0)$, and the fact that $\LL$ is a non-decreasing function assure that for all $t\in[0,\Ts)$ it holds that
\begin{align}\label{eq:2_32}
\norm{ \V(t) }_{\fh,2} & \le \exp\pr[\big]{ t ( \LL_t + \logNorm_{\fh,2}(A) ) } \norm{ \V(0) }_{\fh,2} + \int_0^t \exp\pr[\big]{ (t-w) ( \LL_t + \logNorm_{\fh,2}(A) ) } \norm{C(w)}_{\fh,2} \dx w \nonumber \\
& \le \exp\pr[\big]{ t ( \LL_t + \logNorm_{\fh,2}(A) ) } \int_0^t \norm{C(w)}_{\fh,2} \dx w \dpp
\end{align}
Next, observe that the fact that for all $j \in \N$, $v_1,v_2,\ldots,v_j \in [0,\infty)$ it holds that $[\sum_{k=1}^j v_j]^{\nicefrac{1}{2}} \le \sum_{k=1}^j (v_j)^{\nicefrac{1}{2}}$ guarantees that for all $s\in[0,\Ts)$ it holds that
\begin{equation}\label{eq:2_33}
\begin{split}
\norm{\CC(s)}_{\fh,2} 
& = \left[ \SmallSum_{k=1}^N  \tfrac{1}{2} (h_{k-1} + h_k) \abs[\big]{ \CC_k(s) }^2 \right]^{\!\nicefrac{1}{2}} \\
& \le \left[ \SmallSum_{k=1}^N \tfrac{1}{2} (h_{k-1} + h_k) \right]^{\!\nicefrac{1}{2}} \left[ \max_{k\in\{1,2,\ldots,N\}} \abs{ \CC_k(s) } \right] 
\le \sqrt{2a} \left[ \max_{k\in\{1,2,\ldots,N\}} \abs{ \CC_k(s) } \right] \dpp
\end{split}
\end{equation}
Moreover, note that the triangle inequality shows that for all $s\in[0,\Ts)$, $n\in\{1,2,\ldots,N\}$ it holds that
\begin{align}
\abs{ \CC_n(s) } 
& = \abs*{ - \int_0^{h_{n-1}} \frac{\pr[\big]{ \tfrac{\partial^3}{\partial x^3} u}\lrSpace(t,s)}{h_{n-1}(h_{n-1}+h_n)} (h_{n-1}-s)^2 \dx s +  \int_0^{h_{n}} \frac{\pr[\big]{ \tfrac{\partial^3}{\partial x^3} u}\lrSpace(t,s)}{h_{n}(h_{n-1}+h_n)} (h_{n}-s)^2 \dx s } \nonumber \\
& \le \max_{k\in\{n-1,n\}} \abs*{ \int_0^{h_k} \pr[\big]{ \tfrac{\partial^3}{\partial x^3} u}\lrSpace(t,s) \dx s } \dpp
\end{align}
Combining this, \cref{eq:2_32}, \cref{eq:2_33}, and \cref{lem:A_log} proves that for all $t\in[0,\Ts)$ it holds that
\begin{equation}
\begin{split}
\norm{ \V(t) }_{\fh,2} 
& \le \exp\pr[\big]{ t ( \LL_t + \logNorm_{\fh,2}(A) ) } \int_0^t \sqrt{2a} \left[ \max_{n\in\{1,2,\ldots,N\}} \abs{ C_n(w) } \right] \dx w \\
& \le \sqrt{2a} \exp(\Ts \LL_t) \left[ \int_0^t \max_{n\in\{1,2,\ldots,N\}} \abs*{ \int_0^{h_n} \pr[\big]{ \tfrac{\partial^3}{\partial x^3} u}\lrSpace(w,s) \dx s } \dx w \right] \dpp
\end{split}
\end{equation}
\end{aproof}

We close this subsection with some comments regarding \cref{eq:sd_conv_bd_eq} in \cref{lem:space_conv}.
It is important to note the bound in \cref{eq:sd_conv_bd_eq} could potentially be degenerate if the solution $u$ to \cref{pde2} does not possess sufficient regularity.
We will further estimate the bound in \cref{eq:sd_conv_bd_eq} in \cref{sec:5}; however, it should be clear to readers that if we have that, e.g., $u \in C^{1,3}([0,\Ts),[-a,a])$ then it holds that we can appropriately bound the expression in \cref{eq:sd_conv_bd_eq}.
Moreover, we can obtain similarly useful bounds under the less restrictive assumption that, e.g., $u \in C^{1,2}([0,\Ts),[-a,a])$ with $(\tfrac{\partial^2}{\partial x^2}u) \colon [0,\Ts)\times[-a,a] \to \R$ being absolutely continuous.

\section{The nonlinear operator splitting method}\label{sec:4}

In this section, we introduce the fully-discrete method to be employed in approximating the solution $u$ to \cref{pde2}.
This method is fully implicit and based upon the notion of \emph{nonlinear operator splitting}.
We motivate the general idea behind nonlinear operator splitting in \cref{sec:4_1} below.
These concepts are well-understood in the abstract setting, but there is still much work to be done when implementing such methods for problems exhibiting singularities (e.g., in the case of quenching-combustion PDEs).
Next, in \cref{sec:4_2} below, we introduce the proposed numerical algorithm in \cref{setting3}.
Thereafter, in \cref{sec:4_3} below, we explore the positivity and monotonicity of the sequence of approximations obtained via the proposed nonlinear operator splitting method (cf.\ \cref{eq:fully_discrete} below).

\subsection{General motivation of nonlinear operator splitting methods}\label{sec:4_1}

We now briefly motivate the ideas behind the nonlinear operator splitting method which is implemented herein. 
For more details, we refer interested readers to, e.g., \cite{MR2960375,MR2670117,MR3217072,MR1469677,MR2997899} and the references therein.
Let $T \in (0,\infty)$,
let $(X,\norm{\cdot})$ be a normed $\C$-Banach space, 
let $B(X)$ be the set of all bounded operators on $X$,
let $u_0 \in X$,
let $A \colon D(A) \subseteq X \to X$ be linear and the generator of a strongly-continuous (contraction) semigroup,
let $\semi_t(A) \colon X \to B(X)$, $t\in[0,T]$, be the strongly-continuous semigroup generated by $A$,
let $F \colon D(F) \subseteq X \to X$ be (possibly) nonlinear, 
and let $u \colon [0,T] \to X$ satisfy for all $t\in[0,T]$ that $u(0) = u_0$ and
\begin{equation}\label{motivate0}
\pr[\big]{ \tfrac{{\rm d}}{{\rm d}t} u }\lrSpace(t) = (Au)(t) + (Fu)(t) \dpp
\end{equation}
Then it is a well-known result that (under technical regularity assumptions) for all $t\in[0,T]$ it holds that
\begin{equation}\label{motivate1}
u(t) = \semi_t(A) u_0 + \int_0^t \semi_{t-s}(A) (Fu)(s) \dx s \dpp
\end{equation}
The general idea of nonlinear operator splitting is to approximate \cref{motivate1} via a ``linear kick'' followed by a ``nonlinear kick.''
That is, 
let $T_* \in (0,T]$ and
let $v \colon [0,T_*] \to X$, $w \colon [0,T_*] \to X$ satisfy for all $t\in[0,T_*]$ that $v(0) = u_0$, $w(0) = v(T_*)$,
\begin{equation}\label{motivate3}
\pr[\big]{ \tfrac{{\rm d}}{{\rm d}t} v }\lrSpace(t) = (Av)(t) \dc
\qquad \text{and} \qquad 
\pr[\big]{ \tfrac{{\rm d}}{{\rm d}t} w }\lrSpace(t) =  (Fw)(t) \dpp
\end{equation}
Then (under more technical regularity assumptions) it holds that there exists $\mathfrak{C} \in \R$ such that
\begin{equation}\label{motivate2}
w(t) = \semi_{T_*}(A)u_0 + \int_0^t (Fw)(s) \dx s
\qquad \text{and} \qquad
\norm{ u(T_*) - w(T_*) } \le \mathfrak{C} (T_*)^2 \dpp
\end{equation}

Observe that in \cref{motivate2} it has been assumed that we can solve the sub-problems in \cref{motivate3} \emph{exactly}.
Therefore, we may also consider \emph{approximations} to the sub-problems in \cref{motivate3}. 
The derivation of such approximations is precisely the goal of the current section.
There are many possible methods for approximating semigroups, however, the most common choice is the use of so-called Pad{\'e} or rational approximations
(cf., e.g., \cite{MR537280,MR899800,MR2261021}).
Moreover, based on the results outlined in \cref{sec:quench_props}, it should be clear that we are only interested in approximation methods which preserve the \emph{positivity} and \emph{monotonicity} of the true solution $u$ to \cref{pde2}.
This goal of preserving these qualitative features of the solution to \cref{pde2} is explored further in \cref{sec:4_3}.

\begin{remark}
It is worth noting that we have considered only the one-dimensional problem for exactly the reasons outlined above.
That is, if we had considered, e.g., the two-dimensional quenching-combustion problem on a rectangular domain, then the semidiscretized solution would still have the form \cref{motivate0}, where $A = A_1 + A_2$.
One can then apply classical dimensional splitting to $\semi_{T_*}(A_1+A_2)$ (cf., e.g., \cite{MR2501059,MR2369204,MR1000457,Sheng2}).
The analysis of each of the resulting sub-problems then follows in a fashion similar to the one-dimensional case presented below.
\end{remark}

\subsection{The fully-discretized method}\label{sec:4_2}

In this subsection, we introduce the fully-discrete nonlinear operator splitting algorithm to be employed in the remainder of this article. 
The details of this algorithm are outlined in \cref{setting3} below.
Note that an important feature of the algorithm is the conditions on the nonuniform temporal grid given in \cref{eq:time_step} below.
We close the subsection by providing some elementary results regarding an alternative representation of the numerical algorithm in \cref{eq:fully_discrete} in \cref{lem:mon1a} below.

\begin{setting}\label{setting3}
Assume \cref{setting2},
let $\delta \in (0,\min\{1,\Ts\})$,
$\tau_{-1} , \tau_0,\tau_1,\tau_2,\ldots \in [0,\Ts)$, $t_0 , t_1, t_2 , \ldots \in [0,\infty)$, $\temp{0},\temp{1},\temp{2},\ldots \in \R^N$ satisfy for all $k \in \N_0$, $n \in \{1,2,\ldots,N\}$ that $\tau_{-1} = 1$, $t_0 = 0$, $t_{k+1} = t_k + \tau_k$, $\temp{0}_n = U_n(0)$,
\begin{equation}\label{eq:time_step}
\frac{\tau_k}{\delta} 
= \min_{i\in\{1,2,\ldots,N\}}\pr*{ \min\cu*{ \frac{ 1 - \temp{k+1}_i }{ f(\temp{k+1}_i) } \dc \frac{ 1 }{ f'(\temp{k}_i) } } }
\dc
\end{equation}
and
\begin{equation}\label{eq:fully_discrete}
\temp{k+1} - \temp{k} + (\tau_k)^2 AF(\temp{k+1}) = \tau_k A \temp{k+1} + \tau_k F(\temp{k+1}) \dpp
\end{equation}
\end{setting}

We now clarify the conditions imposed on the sequence $\{\tau_k\}_{k\in\N_0}$ in \cref{eq:time_step} in \cref{setting3}.
Observe that the first condition in \cref{eq:time_step} is only needed to prove that for all $k\in\N_0$, $n \in \{1,2,\ldots,N\}$ it holds that $\temp{k}_n < 1$. 
Again, this condition can be removed under the tacit assumption that for all $k\in\N_0$, $n \in \{1,2,\ldots,N\}$ it holds that $\tau_k$ is sufficiently small to ensure that $\temp{k}_n < 1$ (this is the standard approach in existing numerical approaches to quenching-combustion PDEs).
However, this implicit condition can be updated with a more restrictive explicit condition, if desired (cf., e.g., \cref{eq:upper_sol_bd00} in the proof of \cref{lem:fd_sol_bd} below).
Moreover, the \emph{implicit} nature of the condition on the temporal grid has occurred frequently in the approximation of quenching-combustion PDEs (cf., e.g, the so-called arc-length monitors employed in Sheng \cite{MR2563643}).
Next, note that the second condition in \cref{eq:time_step} is the truly necessary condition.
This condition is the one which allows one to prove that the numerical approximations produced by \cref{eq:fully_discrete} recover the desired monotonicity property of the true solution.
Finally, it should be clear that \cref{eq:time_step} and the fact that $[0,1) \ni s \mapsto \nicefrac{(1-s)}{f(s)} \in [0,\infty)$ is non-increasing assure that for all $k \in \N_0$ it holds that $\tau_{k+1} \le \tau_k$.

While not discussed herein, it is often the case in practice that one will simply allow the initial time steps to be uniform and then adaptation will only commence once the maximum value of the numerical approximation reaches some threshold value (in fact, doing so would not change any of the ensuing results, so long as the step sizes are small enough to prevent \emph{numerical quenching}).

\begin{athm}{lemma}{lem:mon1a}
Assume \cref{setting3}.
Then 
\begin{enumerate}[label=(\roman*)]
\item\label{lem:mon1a_i1} for all $k \in \N_0$ it holds that $\idMatrix_N - \tau_kA$ is nonsingular and
\item\label{lem:mon1a_i2} for all $k \in \N_0$ it holds that 
$\temp{k+1} = (\idMatrix_N - \tau_k A)^{-1} \temp{k} + \tau_k F(\temp{k+1})$
\end{enumerate}
(cf.\ \cref{def:identityMatrix}).
\end{athm}

\begin{aproof}
First, note that \cref{eq:mat_A} ensures that for all $k\in\N_0$, $i,j,n\in\{1,2,\allowbreak\ldots,\allowbreak N\}$, $l\in\{1,2,\ldots,N-1\}$ with $\abs{i-j} \in \{2,3,\ldots,N-1\}$ it holds that
\begin{equation}\label{lem:mon1a_1}
\pr[\big]{ \idMatrix_N - \tau_k A}_{i,j} = 0
\dc \qquad
\pr[\big]{ \idMatrix_N - \tau_k A}_{n,n} = 1 + \tfrac{2\tau_k}{h_{n-1}h_n} \dc
\end{equation}
\begin{equation}\label{lem:mon1a_2}
\pr[\big]{ \idMatrix_N - \tau_k A}_{l+1,l} =  -\tfrac{2\tau_k}{h_{l}(h_l + h_{l+1})}
\dc
\qquad \text{and} \qquad
\pr[\big]{ \idMatrix_N - \tau_k A}_{l,l+1} = -\tfrac{2\tau_k}{h_{l}(h_{l-1}+h_l)} \dpp
\end{equation}
Combining \cref{lem:mon1a_1,lem:mon1a_2} hence assures that for all $k\in\N_0$, $n\in\{1,2,\ldots,N\}$ it holds that
\begin{equation}
\SmallSum_{\substack{j=1\\j\neq n}}^N \abs[\big]{ \pr[\big]{ \idMatrix_N - \tau_k A}_{n,j} }
= \tfrac{2\tau_k}{h_{n-1}(h_{n-1} + h_{n})} + \tfrac{2\tau_k}{h_{n}(h_{n-1}+h_n)} 
< 
\abs[\big]{ \pr[\big]{ \idMatrix_N - \tau_k A}_{n,n} } \dpp
\end{equation}
This and, e.g., Horn and Johnson \cite[Theorem 6.1.10]{Horn} establish \cref{lem:mon1a_i1}.
Next, observe that combining \cref{lem:mon1a_i1} and \cref{eq:fully_discrete} proves \cref{lem:mon1a_i2}.
\end{aproof}

\subsection{Positivity and monotonicity of the operator splitting method}\label{sec:4_3}

In this subsection, we study the componentwise positivity and monotonicity of the approximation to the solution $u$ to \cref{pde2} obtained via \cref{eq:fully_discrete}.
First, we prove the positivity of the numerical approximations obtained via \cref{eq:fully_discrete} in \cref{lem:sol_sd_pos}.
The proof of this result in turn is based on \cref{lem:pos_mat_vec}.
\cref{lem:pos_mat_vec} demonstrates that our Pad{\'e} approximation of the underlying matrix exponential is a \emph{monotone} matrix (cf.\ \cref{def:mat_exp}).
Next, we prove \cref{lem:fd_sol_bd} which demonstrates that the components of the numerical approximations obtained via \cref{eq:fully_discrete} remain bounded componentwise by unity for all $k\in\N_0$.
The proof of \cref{lem:fd_sol_bd} in turn depends on the elementary positivity result in \cref{prop:mat_bd}.
Finally, we prove \cref{lem:mon2}, which is the main result of this subsection.
\cref{lem:mon2} demonstrates that the sequence of approximations obtained via \cref{eq:fully_discrete} is actually a sequence which is monotonically increasing (componentwise) with respect to $k\in\N_0$.

\begin{athm}{lemma}{lem:pos_mat_vec}
Assume \cref{setting3} and let $w \in \R^N$ be nonnegative (cf.\ \cref{def:wedge}).
Then it holds for all $t\in[0,\infty)$ that $(\idMatrix_N - t A)^{-1} w$ is nonnegative (cf.\ \cref{def:identityMatrix}).
\end{athm}

\begin{aproof}
Throughout this proof let $\beta = \nicefrac{2}{\min_{n\in\{1,2,\ldots,N\}} h_{n-1}h_n}$ and let $B = A + \beta \idMatrix_N$.
First, observe that for all $t\in[0,\infty)$, $i,j\in\{1,2,\ldots,N\}$ with $i\neq j$ it holds that
\begin{equation}\label{lem:pos_mat_vec_eq1}
\pr[\big]{ \idMatrix_N - t A }_{i,j} \le 0 \dpp
\end{equation}
Next, note that for all $t\in[0,\infty)$ it holds that
\begin{equation}\label{lem:pos_mat_vec_eq2}
\idMatrix_N - t A = \idMatrix_N - t (B - \beta \idMatrix_N) = (1+\beta t) \idMatrix_N - t B \dpp
\end{equation}
Combining \cref{lem:pos_mat_vec_eq1}, \cref{lem:pos_mat_vec_eq2}, the fact that for all $i,j\in\{1,2,\ldots,N\}$ it holds that $B_{i,j} \ge 0$, and, e.g., 
Plemmons \cite[Condition F$_{15}$]{MR444681} guarantees that for all $t\in[0,\infty)$, $w\in\R^N$ with $w \ge 0$ it holds that $(\idMatrix_N - t A)^{-1} w$ is nonnegative.
\end{aproof}

\begin{athm}{lemma}{lem:sol_sd_pos}
Assume \cref{setting3}.
Then it holds for all $k\in\N_0$ that $\temp{k}$ is nonnegative (cf.\ \cref{def:wedge}).
\end{athm}

\begin{aproof}
Throughout this proof let $\0 \in \R^N$ satisfy for all $n\in\{1,2,\ldots,N\}$ that $\0_n = 0$ and let $W^{(k)} \colon [0,\tau_{k-1}] \allowbreak \to \R^N$, $k\in\N$, satisfy for all $k\in\N$, $t\in[0,\tau_{k-1}]$ that
\begin{equation}
W^{(k)}(t) = (\idMatrix_N - tA)^{-1} \temp{k-1} + t F(W^{(k)}(t))
\end{equation}
(cf.\ \cref{def:identityMatrix}).
Note that the assumption that $f\in C([0,1),\R)$ and \cref{lem:mon1a_i2} in \cref{lem:mon1a} ensure that
\begin{enumerate}[label=(\Roman*)]
\item\label{proof_i1} it holds for all $k\in\N$ that $W^{(k)} \in C([0,\tau_{k-1}],\R^N)$ and
\item\label{proof_i2} it holds for all $k\in\N$ that $W^{(k)}(\tau_{k-1}) = \temp{k}$.
\end{enumerate}
Next, observe that the assumption that for all $x\in[-a,a]$ it holds that $u_0(x) \in [0,1)$ assures that $\temp{0}$ is nonnegative (cf.\ \cref{def:wedge}).
Moreover, we claim that for all $k\in\N$, $t\in[0,\tau_{k-1}]$ it holds that
\begin{equation}\label{pos_induct}
W^{(k)}(t) > 0 \dpp
\end{equation}
We prove \cref{pos_induct} by induction on $k\in\N$. 
For the sake of a contradiction, we claim that there exists $\tau_* \in [0,\tau_0]$ such that for all $t \in [0,\tau_*)$ it holds that $W^{(1)}(t)$ is positive and $W^{(1)}(\tau_*) = \0$.
Note that this, the fact that $\temp{0}$ is nonnegative, \cref{lem:pos_mat_vec}, \cref{proof_i1}, and the assumption that $f(0) > 0$ guarantee that
\begin{equation}
\0 = W^{(1)}(\tau_*) = (\idMatrix_N - \tau_*A)^{-1} \temp{0} + \tau_* F(W^{(1)}(\tau_*)) > 0 \dpp
\end{equation}
This is a contradiction to the claim that $W^{(1)}(\tau_*) = \0$. 
This establishes \cref{pos_induct} in the base case $k=1$.
For the induction step $\N \ni (k-1) \induct k \in \{2,3,4,\ldots\}$, let $k\in\{2,3,4,\ldots\}$ and assume that for every $i\in\{1,2,\ldots,k-1\}$, $t\in[0,\tau_{i-1}]$ it holds that $W^{(i)}(t)$ is positive.
For the sake of a contradiction, we claim that there exists $\tau_* \in [0,\tau_{k-1}]$ such that for all $t\in[0,\tau_*)$ it holds that $W^{(k)}(t)$ is positive and $W^{(k)}(\tau_*) = \0$.
Observe that this, the induction hypothesis, \cref{lem:pos_mat_vec}, \cref{proof_i1}, and the assumption that $f(0) > 0$ demonstrate that
\begin{equation}
\0 = W^{(k)}(\tau_*) = (\idMatrix_N - \tau_*A)^{-1} \temp{k-1} + \tau_* F(W^{(k)}(\tau_*)) > 0 \dpp
\end{equation}
This is a contradiction to the assumption that $W^{(k)}(\tau_*) = \0$. 
This and induction thus establish \cref{pos_induct}.
In addition, note that combining \cref{proof_i2} and the fact that $\temp{0}$ is nonnegative proves that for all $k\in\N_0$ it holds that $\temp{k}$ is nonnegative.
\end{aproof}

\begin{remark}
It should be clear from the proof of \cref{lem:sol_sd_pos} that for all $k\in\N$ it holds that $\temp{k}$ is positive (cf.\ \cref{def:wedge}). 
Moreover, if $\temp{0}$ is positive, then for all $k\in\N_0$ it holds that $\temp{k}$ is positive.
\end{remark}

\begin{athm}{proposition}{prop:mat_bd}
Assume \cref{setting3}, let $\fc \in [0,\infty)$, and let $\fx \in \R^N$ satisfy that $\fx = (\fc, \fc, \ldots,\fc)^*$ (cf.\ \cref{def:Transpose}).
Then it holds for all $t\in[0,\infty)$ that
\begin{equation}\label{prop:mat_bd_eq}
(\idMatrix_N - t A)^{-1} \fx \ge \fx
\end{equation}
(cf.\ \cref{def:identityMatrix,def:wedge}).
\end{athm}

\begin{aproof}
Note that \cref{prop:mat_bd_eq} follows immediately from \cref{eq:mat_A}.
\end{aproof}

\begin{athm}{lemma}{lem:fd_sol_bd}
Assume \cref{setting3}.
Then it holds for all $k\in\N_0$ that
\begin{equation}\label{eq:upper_sol_bd_goal}
\max_{n\in\{1,2,\ldots,N\}} \temp{k}_n \le 1 - (1+\delta)^{-k} (1 - \temp{0}_n ) < 1 \dpp
\end{equation}
\end{athm}

\begin{aproof}
Throughout this proof let $\fx \in \R^N$ satisfy that $\fx = (1, 1, \ldots,1)^*$ (cf.\ \cref{def:Transpose}).
We claim that for all $k\in\N_0$ it holds that
\begin{equation}\label{eq:upper_sol_bd}
\temp{k+1} - \fx = (1+\delta)^{-(k+1)} \br*{ \prod_{j=0}^k (\idMatrix_N - \tau_j A)^{-1} } (\temp{0} - \fx) 
\end{equation}
(cf.\ \cref{def:identityMatrix}).
We prove \cref{eq:upper_sol_bd} by induction on $k\in\N_0$.
Note that \cref{lem:mon1a_i2} in \cref{lem:mon1a} assures that for all $k\in\N_0$ it holds that
\begin{equation}
\begin{split}
\temp{k+1} - \fx 
& = \br[\Big]{ (\idMatrix_N - \tau_{k} A)^{-1} \temp{k} + \tau_{k} F(\temp{k+1}) } - \fx \\
& = (\idMatrix_N - \tau_{k} A)^{-1} \br[\big]{ \temp{k} - (\idMatrix_N -\tau_{k} A) \fx } + \tau_{k} F(\temp{k+1}) 
\dpp
\end{split}
\end{equation}
Combining this, \cref{prop:mat_bd}, and \cref{eq:time_step} demonstrates that for all $k\in\N_0$ it holds that
\begin{equation}\label{eq:upper_sol_bd0}
\begin{split}
\temp{k+1} - \fx 
& \le (\idMatrix_N - \tau_{k} A)^{-1} \br[]{ \temp{k} - \fx } + \tau_{k} F(\temp{k+1}) \\
& \le (\idMatrix_N - \tau_{k} A)^{-1} \br[]{ \temp{k} - \fx } - \delta ( \temp{k+1} - \fx ) 
\end{split}
\end{equation}
(cf.\ \cref{def:wedge}).
This implies that for all $k\in\N_0$ it holds that
\begin{equation}\label{eq:upper_sol_bd00}
\temp{k+1} - \fx 
\le (1+\delta)^{-1} (\idMatrix_N - \tau_{k} A)^{-1} \br[]{ \temp{k} - \fx }
\dpp
\end{equation}
Observe that \cref{eq:upper_sol_bd00} and the assumption that $\delta \in (0,1)$ hence establish \cref{eq:upper_sol_bd} in the base case $k=0$.
For the induction step $\N_0 \ni (k-1) \induct k \in \N$, let $k\in\N$ and assume that for every $i\in\{0,1,\ldots,k-1\}$ that \cref{eq:upper_sol_bd} holds.
Combining the induction hypothesis and \cref{eq:upper_sol_bd00} therefore demonstrates that
\begin{equation}\label{eq:upper_sol_bd1}
\begin{split}
\temp{k+1} - \fx 
& \le (1+\delta)^{-1} (\idMatrix_N - \tau_{k} A)^{-1} \br[]{ \temp{k} - \fx } \\
& \le (1+\delta)^{-1} \pr*{ (1+\delta)^{-k} \br*{ \prod_{j=0}^{k-1} (\idMatrix_N - \tau_{k} A)^{-1} (\idMatrix_N - \tau_j A)^{-1} } (\temp{0} - \fx) } \\
& = (1+\delta)^{-(k+1)} \br*{ \prod_{j=0}^k (\idMatrix_N - \tau_j A)^{-1} } (\temp{0} - \fx)
\dpp
\end{split}
\end{equation}
Induction hence establishes \cref{eq:upper_sol_bd}.
Next, note that combining \cref{eq:upper_sol_bd1}, the fact that for all $n\in\{1,2,\ldots,N\}$ it holds that $\temp{0}_n \in [0,1)$, and \cref{prop:mat_bd} 
(applied for every $j\in\{0,1,\ldots,k\}$ with $t \with \tau_j$, $\idMatrix_N \with \idMatrix_N$, $A \with A$, $\fx \with \pr[]{ \br[]{ {\textstyle\min_{n\in\{1,2,\ldots,N\}} \temp{0}_n} } \idMatrix_N  - \fx }$ in the notation of \cref{prop:mat_bd})
guarantees that for all $k\in\N_0$ it holds that
\begin{equation}
\begin{split}
\temp{k+1} - \fx 
& \le (1+\delta)^{-(k+1)} \br*{ \prod_{j=0}^k (\idMatrix_N - \tau_j A)^{-1} } \pr[\Big]{ \br[\big]{ {\textstyle\min_{n\in\{1,2,\ldots,N\}} \temp{0}_n} } \idMatrix_N  - \fx } \\
& \le (1+\delta)^{-(k+1)} \pr[\Big]{ \br[\big]{ {\textstyle\min_{n\in\{1,2,\ldots,N\}} \temp{0}_n} } \idMatrix_N  - \fx }
\dpp
\end{split}
\end{equation}
Combining this, the fact that $\delta\in[0,1)$, and the fact that for all $n\in\{1,2,\ldots,N\}$ it holds that $\temp{0}_n \in [0,1)$
proves \cref{eq:upper_sol_bd_goal}.
\end{aproof}

\begin{athm}{lemma}{lem:mon2}
Assume \cref{setting3} and assume that $A\temp{0} + F(\temp{0}) - \tau_0 A F(\temp{0})$ is nonnegative (cf.\ \cref{def:wedge}).
Then it holds for all $k \in \N_0$
that $\temp{k+1} - \temp{k}$ is nonnegative.
\end{athm}

\begin{aproof}
Throughout this proof 
let $\w{0} , \w{1} , \w{2} , \ldots \in\R^N$ satisfy for all $k\in\N_0$ that $\w{k} = (\tau_k)^{-1} [\temp{k+1} - \temp{k}]$
and
for every $G \colon \R^N \to \R^N$ let $J_G \colon \R^N \to \R^{N\times N}$ satisfy for all $w \in \R^N$, $i,j\in\{1,2,\ldots,N\}$ that
\begin{equation}\label{eq:jac_def}
(J_G(w))_{i,j} = \pr[\big]{ \tfrac{\partial}{\partial x_j} G_i } \lrSpace (w) \dpp
\end{equation}
We claim that for all $k\in\N_0$ it holds that
\begin{equation}\label{eq:mon_claim}
\w{k} \ge 0 \dpp
\end{equation}
We now prove \cref{eq:mon_claim} by induction on $k\in\N_0$.
For the base case $k =0$, note that \cref{lem:mon1a_i2} in \cref{lem:mon1a} guarantees that
\begin{equation}\label{eq:mon_claim_1}
\tau_0 \w{0} 
= (\idMatrix_N - \tau_0 A)^{-1} \temp{0} + \tau_0 F(\temp{1}) - \temp{0} 
= \tau_0 (\idMatrix_N - \tau_0 A)^{-1} A \temp{0} + \tau_0 F(\temp{1})
\end{equation}
(cf.\ \cref{def:identityMatrix}).
Next, observe that the hypothesis that $f \in C^1([0,1),\R)$, \cref{f_cond}, and \cref{eq:jac_def}
assure that for all $k\in\N_0$ 
it holds that
\begin{equation}\label{eq:F_jac}
F(\temp{k+1}) 
\ge F(\temp{k}) + J_F(\temp{k}) (\temp{k+1} - \temp{k} ) 
= F(\temp{k}) + \tau_k J_F(\temp{k}) \w{k}
\dpp
\end{equation}
Combining this and \cref{eq:mon_claim_1} ensure that
\begin{equation}
\begin{split}
\w{0} 
& \ge (\idMatrix_N - \tau_0 A)^{-1} A \temp{0} + F(\temp{0}) + \tau_0 J_F(\temp{0}) \w{0} \\
& = (\idMatrix_N - \tau_0 A)^{-1} \br[\big]{ A\temp{0} + F(\temp{0}) - \tau_0 A F(\temp{0}) } + \tau_0 J_F(\temp{0}) \w{0}
\dpp
\end{split}
\end{equation}
This, the assumption that $A\temp{0} + F(\temp{0}) - \tau_0 A F(\temp{0})$ is nonnegative, the fact that for all $x\in[0,1)$ it holds that $f'(x) > 0$, \cref{eq:time_step}, \cref{lem:pos_mat_vec}, \cref{lem:sol_sd_pos}, and \cref{lem:fd_sol_bd} guarantee that
\begin{equation}
\w{0} \ge \pr[\big]{ \idMatrix_N - \tau_0 J_F(\temp{0}) }^{-1} (\idMatrix_N - \tau_0 A)^{-1} \br[\big]{ A\temp{0} + F(\temp{0}) - \tau_0 A F(\temp{0}) } \ge 0 \dpp
\end{equation}
This establishes \cref{eq:mon_claim} in the base case $k=0$.
For the induction step $\N_0 \ni (k-1) \induct k \in \N$, let $k\in\N$ and assume that for all $i \in \{0,1,\ldots,k-1\}$ it holds that $\w{i} \ge 0$.
Note that \cref{eq:fully_discrete} ensures that 
\begin{equation}\label{eq:4_29}
\w{k} - \w{k-1} - \tau_{k} A \w{k}
= (\idMatrix_N - \tau_{k} A) F(\temp{k+1}) - (\idMatrix_N - \tau_{k-1} A) F(\temp{k}) 
\dpp
\end{equation}
Combining this, \cref{eq:time_step}, and \cref{eq:F_jac} shows that
\begin{align}
& (\idMatrix_N - \tau_{k} A) F(\temp{k+1}) - (\idMatrix_N - \tau_{k-1} A) F(\temp{k}) \nonumber \\
& = \br[\big]{ F(\temp{k+1}) - F(\temp{k}) } - A \br[\big]{ \tau_{k} F(\temp{k+1}) - \tau_{k-1} F(\temp{k}) } \\
& \ge \tau_{k} J_F(\temp{k}) \w{k} - \tau_{k} A \br[\big]{ F(\temp{k+1}) - \tfrac{\tau_{k-1}}{\tau_{k}} F(\temp{k}) } 
\ge \tau_{k} J_F(\temp{k}) \w{k} - \tau_{k} A J_F(\temp{k}) \w{k} 
\dpp \nonumber
\end{align}
This and \cref{eq:4_29} demonstrate that 
\begin{equation}
\begin{split}
\w{k-1} & \le \w{k} - \tau_{k} A \w{k} - \tau_{k} J_F(\temp{k}) \w{k} + \tau_{k} A J_F(\temp{k}) \w{k} \\
& = (\idMatrix_N - \tau_k A) \pr[\big]{ \idMatrix_N - \tau_k J_F(\temp{k}) } \w{k} \dpp
\end{split}
\end{equation}
Combining this, the fact that for all $x\in[0,1)$ it holds that $f'(x) > 0$, \cref{eq:time_step}, \cref{lem:pos_mat_vec}, \cref{lem:sol_sd_pos}, \cref{lem:fd_sol_bd}, and the inductive hypothesis prove that
\begin{equation}
\w{k} \ge \pr[\big]{ \idMatrix_N - \tau_k J_F(\temp{k}) }^{-1} (\idMatrix_N - \tau_k A)^{-1} \w{k-1} \ge 0 \dpp
\end{equation}
Induction hence establishes \cref{eq:mon_claim}.
Moreover, observe that \cref{eq:mon_claim} and the assumption that for all $j\in\N_0$ it holds that $\tau_j \in [0,\Ts)$ yield that for all $k\in\N_0$ it holds that
\begin{equation}
0 \le \tau_k \w{k} = \temp{k+1} - \temp{k} \dpp
\end{equation}
\end{aproof}

\begin{remark}
Note that \cref{lem:mon2} demonstrates that under the assumptions outlined in \cref{setting3}, it holds that $\{\temp{k}\}_{k\in\N_0}$ is a monotonically increasing sequence (where ``monotonically increasing'' is defined in the componentwise fashion outlined in \cref{def:wedge}).
\end{remark}

\section{Convergence analysis of the nonlinear operator splitting method}\label{sec:5}

In this section, we carry out the necessary space-time convergence analysis of the proposed implicit nonlinear operator splitting method (cf.\ \cref{eq:fully_discrete} above).
A primary contribution of the ensuing analysis is the emphasis on exploring the degree to which the quenching singularity may impede the expected rate of convergence.
Such impedance has been observed in numerical experiments, but has yet to be quantified precisely through rigorous analysis (cf.\ Padgett and Sheng \cite[Section 6]{MR3759137}). 

First, in \cref{lem:pre_converge} in \cref{sec:exp_bd} below, we develop bounds for certain terms which arise in the analysis of the error between the semidiscretized and fully-discretized numerical solutions.
The proof of \cref{lem:pre_converge} depends on
the discrete weighted $2$-norm bounds of the matrix exponential and its first-order subdiagonal Pad{\'e} approximation developed in \cref{lem:exp_bd}.
Next, in \cref{sec:5_2} below, we develop a representation of the error between the semidiscrete and fully-discrete numerical solution.
The main result of this subsection is the global error result presented in \cref{lem:converge2}.
The temporal global error bound developed in \cref{lem:converge2} is a direct consequence of the temporal local error bound developed in \cref{lem:converge1}.
Finally, in \cref{th:final} in \cref{sec:5_3} below, we present the main results of this article.
In \cref{th:final} below, we present a general result regarding the space-time convergence rates of the proposed implicit nonlinear operator splitting method.
These results are then used to develop \cref{cor:final} below, which interprets our results in \cref{th:final} in the setting originally studied by Hideo Kawarada \cite{Kawa}.

\subsection{The exponential matrix and its approximation}\label{sec:exp_bd}

In this subsection, we develop preliminary results needed for the convergence analysis carried out in \cref{sec:5_2} below.
Note that \cref{lem:exp_bd} is an extension of well-known results regarding bounds of the matrix exponential and its first-order subdiagonal Pad{\'e} approximation via the logarithmic norm.
Interested readers can see, e.g., \cite{jones2020intrinsic,
MR2265579,MR408227} for bounds related to the matrix exponential
and, e.g., \cite{MR537280,MR1201070,MR899800,MR2261021,MR3127028} for a discussion of rational approximations of semigroups.
The results in \cref{lem:exp_bd} are then used in the proofs of \cref{lem:pre_converge} and \cref{lem:converge1} below.

\begin{athm}{lemma}{lem:exp_bd}
Let $d \in \N$, $\fh_0,\fh_1,\ldots,\fh_d \in (0,\infty)$, $\fh = (\fh_0,\fh_1,\ldots,\fh_d)\in \R^{d+1}$. Then 
\begin{enumerate}[label=(\roman*)]
\item\label{lem:exp_bd_i1} it holds for all $v\in\R^d$, $B\in\R^{d\times d}$ that
$\norm{\exp(B) v}_{\fh,2} \le \exp(\logNorm_{\fh,2}(B)) \norm{v}_{\fh,2}$
and
\item\label{lem:exp_bd_i2} it holds for all $v\in\R^d$, $B\in\R^{d\times d}$ with $\logNorm_{\fh,2}(B) \in (-\infty,1)$ that
$\norm{(\idMatrix_d - B)^{-1}v}_{\fh,2} \le (1-\logNorm_{\fh,2}(B))^{-1} \norm{v}_{\fh,2}$
\end{enumerate}
(cf.\ \cref{def:log,def:mat_exp,def:identityMatrix,def:euclid_norm}).
\end{athm}

\begin{aproof}
Throughout this proof 
let $H \in \R^{N\times N}$ be the matrix which satisfies for all $i,j\in\{1,2,\ldots,N\}$ with $i\neq j$ that
$H_{i,i} = \sqrt{\nicefrac{\fh_{k-1}+\fh_k}{2}}$ and $H_{i,j} = 0$,
let $\vt{\cdot,\cdot} \colon \R^N \times \R^N \to \R$ be the function which satisfies for all $x,y \in \R^d$ that $\vt{x,y} = \sum_{k=1}^d x_ky_k$,
without loss of generality let $v \in \R^d$, $B\in\R^{d\times d}$ with $\logNorm_{\fh,2}(B) \in (-\infty,1)$ and $\norm{v}_{\fh,2} \neq 0$,
and let $w \in \R^d$ satisfy that
\begin{equation}\label{eq:log_rat1}
(\idMatrix_d - B) v = w 
\end{equation}
(cf.\ \cref{def:log,def:identityMatrix}).
Note that, e.g., Jones et al.\ \cite[Lemma 2.8]{jones2020intrinsic} establishes \cref{lem:exp_bd_i1}.
Next, observe that \cref{lem:log_rep} and the assumption that $\logNorm_{\fh,2}(B) \in (-\infty,1)$ assure that $\idMatrix_d - B$ is nonsingular.
In addition, note that \cref{eq:log_rat1} guarantees that
\begin{equation}
v - w = Bv \dpp
\end{equation}
This, \cref{eq:inner_to_norm}, \cref{lem:log_rep}, and the Cauchy-Schwartz inequality demonstrate that
\begin{equation}
\begin{split}
\norm{v}_{\fh,2}^2 - \norm{w}_{\fh,2}^2 
& = \vt{H(v+w),H(v-w)} 
= \vt{H(v+w),HBv} \\
& \le \vt{Hw,HBv} + \logNorm_{\fh,2}(B) \norm{v}_{\fh,2}^2 
= \vt{Hw,H(v-w)} + \logNorm_{\fh,2}(B) \norm{v}_{\fh,2}^2 \\
& = \vt{Hw,Hv} - \vt{Hw,Hw} + \logNorm_{\fh,2}(B) \norm{v}_{\fh,2}^2 \\
& \le \norm{w}_{\fh,2} \norm{v}_{\fh,2} - \norm{w}_{\fh,2}^2 + \logNorm_{\fh,2}(B) \norm{v}_{\fh,2}^2 \dpp
\end{split}
\end{equation}
This implies that
\begin{equation}
\pr[\big]{ 1 - \logNorm_{\fh,2}(B) } \norm{v}_{\fh,2} \le \norm{w}_{\fh,2} \dpp
\end{equation}
Combining this, \cref{eq:log_rat1}, the assumption that $\logNorm_{\fh,2}(B) \in (-\infty,1)$, and the fact that $\idMatrix_d - B$ is nonsingular proves that
\begin{equation}
\norm[\big]{(\idMatrix_d - B)^{-1}v}_{\fh,2} \le (1-\logNorm_{\fh,2}(B))^{-1} \norm{v}_{\fh,2} \dpp
\end{equation}
This establishes \cref{lem:exp_bd_i2}.
\end{aproof}

\begin{athm}{lemma}{lem:pre_converge}
Assume \cref{setting3} and let $\fh \in \R^{N+1}$ satisfy $\fh = (h_0,h_1, \ldots, h_N)$.
Then 
\begin{enumerate}[label=(\roman*)]
\item\label{lem:pre_converge_i2a} it holds for all $k\in\N_0$, $s\in[0,\Ts)$ that
\begin{equation}
\norm[\Big]{ \br[\Big]{ \exp\pr[\big]{ (\tau_k-s) A } - \idMatrix_N } F(\temp{k+1}) }_{\fh,2}
\le \tau_k \norm{ AF(\temp{k+1}) }_{\fh,2} 
\end{equation}
and
\item\label{lem:pre_converge_i2} it holds for all $k\in\N_0$ that
\begin{equation}
\norm[\big]{ \pr[\big]{ \exp(\tau_k A) - (\idMatrix_N - \tau_k A)^{-1} } \temp{k} }_{\fh,2}
\le 
2 (\tau_k)^2 \norm{ A^2 \temp{k} }_{\fh,2}
\end{equation}
\end{enumerate}
(cf.\ \cref{def:euclid_norm,def:identityMatrix,def:mat_exp}).
\end{athm}

\begin{aproof}
First, note that, e.g., Fetahu \cite[Theorem 2.12]{fetahu2014semigroups} (applied for every $k\in\N_0$, $s\in[0,\Ts)$ with $A \with A$, $x \with F(\temp{k+1})$, $T(s) \with \exp((\tau_k-s)A)$ in the notation of Fetahu \cite[Theorem 2.12]{fetahu2014semigroups})
ensures that for all $k\in\N_0$, $s\in[0,\Ts)$ it holds that
\begin{equation}
\br[\Big]{ \exp\pr[\big]{ (\tau_k-s) A } - \idMatrix_N } F(\temp{k+1})
= - \int_0^{\tau_k} \exp\pr[\big]{ (\tau_k - \tau) A} A F(\temp{k+1}) \dx \tau
\end{equation}
(cf.\ \cref{def:identityMatrix,def:mat_exp}).
This, Jensen's inequality, \cref{lem:exp_bd_i1} in \cref{lem:exp_bd}, and \cref{lem:A_log} assure that for all $k\in\N_0$, $s\in[0,\Ts)$ it holds that
\begin{align}
& \norm[\Big]{ \br[\Big]{ \exp\pr[\big]{ (\tau_k-s) A } - \idMatrix_N } F(\temp{k+1}) }_{\fh,2}
\le \int_0^{\tau_k} \norm[\big]{ \exp\pr[\big]{ (\tau_k - \tau) A} }_{\fh,2} \norm{ A F(\temp{k+1}) }_{\fh,2} \dx \tau \nonumber \\
& \quad \le \int_0^{\tau_k} \exp\pr[\big]{ (\tau_k - \tau) \logNorm_{\fh,2}(A) }  \norm{ A F(\temp{k+1}) }_{\fh,2} \dx \tau
\le \tau_k \norm{ A F(\temp{k+1}) }_{\fh,2} 
\end{align}
(cf.\ \cref{def:log}).
This establishes \cref{lem:pre_converge_i2a}.
Next, observe that \cref{lem:mon1a_i1} in \cref{lem:mon1a} guarantees that for all $k\in\N_0$ it holds that
\begin{equation}
\pr[\big]{ \exp(\tau_k A) - (\idMatrix_N - \tau_k A)^{-1} } \temp{k}
= (\idMatrix_N - \tau_k A)^{-1} \br[\big]{ (\idMatrix_N - \tau_k A) \exp(\tau_k A) - \idMatrix_N } \temp{k} \dpp
\end{equation}
This and, e.g., Fetahu \cite[Theorem 2.12]{fetahu2014semigroups} (applied for every $k\in\N_0$ with $A \with A$, $x \with \temp{k}$, $t \with \tau_k$, $T(t) \with \exp(\tau_k A)$ in the notation of Fetahu \cite[Theorem 2.12]{fetahu2014semigroups})
yield that for all $k\in\N_0$ it holds that
\begin{equation}\label{eq:pre_converge1}
\begin{split}
& \pr[\big]{ \exp(\tau_k A) - (\idMatrix_N - \tau_k A)^{-1} } \temp{k} \\
& \quad = (\idMatrix_N - \tau_k A)^{-1} \br[\Bigg]{ (\idMatrix_N - \tau_k A) \br[\bigg]{ \idMatrix_N + \int_0^{\tau_k} \exp(sA) A \dx s } - \idMatrix_N } \temp{k} \\
& \quad = (\idMatrix_N - \tau_k A)^{-1} \br[\bigg]{ (\idMatrix_N - \tau_k A) \int_0^{\tau_k} \exp(sA) A \dx s - \tau_k A } \temp{k} \dpp
\end{split}
\end{equation}
Combining \cref{eq:pre_converge1} and integration by parts hence shows that for all $k\in\N_0$ it holds that
\begin{equation}
\begin{split}
& \pr[\big]{ \exp(\tau_k A) - (\idMatrix_N - \tau_k A)^{-1} } \temp{k} \\
& \quad = (\idMatrix_N - \tau_k A)^{-1} \br[\bigg]{ \int_0^{\tau_k} \br[\big]{ \exp(sA) - \tau_k \exp(sA) A - \idMatrix_N } \dx s } A \temp{k} \\
& \quad = (\idMatrix_N - \tau_k A)^{-1} \br[\Bigg]{ \int_0^{\tau_k} \br[\bigg]{ \int_0^s \exp(z A) \dx z - \tau_k \exp(sA) } \dx s } A^2 \temp{k} \dpp
\end{split}
\end{equation}
This, Jensen's inequality, the triangle inequality, \cref{lem:A_log}, and \cref{lem:exp_bd} prove that for all $k \in \N_0$ it holds that
\begin{align}
& \norm[\Big]{ \pr[\big]{ \exp(\tau_k A) - (\idMatrix_N - \tau_k A)^{-1} } \temp{k} }_{\fh,2} \\
& \quad \le \pr[\big]{ 1 - \tau_k \logNorm_{\fh,2}(A) }^{-1} \br[\Bigg]{ \int_0^{\tau_k} \br[\bigg]{ \int_0^s \exp\pr[\big]{ z \logNorm_{\fh,2}(A) } \dx z + \tau_k \exp\pr[\big]{ s \logNorm_{\fh,2}(A) } } \dx s } \norm{ A^2 \temp{k} }_{\fh,2} \nonumber \\
& \quad \le \br[\Bigg]{ \int_0^{\tau_k} \br[\bigg]{ \int_0^s 1 \dx z + \tau_k } \dx s } \norm{ A^2 \temp{k} }_{\fh,2}
\le 2 (\tau_k)^2 \norm{ A^2 \temp{k} }_{\fh,2} \dpp
\nonumber
\end{align}
This establishes \cref{lem:pre_converge_i2}.
\end{aproof}

\subsection{Error representation of the fully-discretized method}
\label{sec:5_2}

In this subsection, we develop the temporal error representation of the proposed fully-discretized numerical solution.
In \cref{lem:converge1} below, we provide a local error bound in the discrete weighted $2$-norm for the difference between the semidiscrete and fully-discrete numerical approximations of quenching-combustion PDEs.
As can be seen in \cref{eplicit_term} in \cref{lem:converge1} below, the impact of the quenching singularity is incorporated into the error bound via the functions $F \colon \R^N \to \R^N$ and $\LL \colon [0,\Ts) \to [0,\infty)$.
We develop the analogous global error bound in the discrete weighted $2$-norm for the difference between the semidiscrete and fully-discrete numerical approximations of quenching-combustion PDEs in \cref{lem:converge2} below.
Again, the impact of the the quenching singularity is incorporated via the functions $F \colon \R^N \to \R^N$ and $\LL \colon [0,\Ts) \to [0,\infty)$
(cf.\ \cref{eq:cor_converge1} in \cref{lem:converge2} below).

\begin{athm}{lemma}{lem:converge1}
Assume \cref{setting3}, let $\fh \in \R^{N+1}$ satisfy $\fh = (h_0,h_1, \ldots, h_N)$, and let $\LL \colon [0,\Ts) \to [0,\infty)$ satisfy for all $t\in[0,\Ts)$ that $\LL_t = \max_{n\in\{1,2,\ldots,N\}} L_{U_n(t),\temp{\min\{k\in\N_0 \colon t \le t_k\}}}$.
Then for all $k \in \N_0$ with $t_{k+1} \in [0,\Ts)$ it holds that
\begin{equation}\label{eplicit_term}
\begin{split}
& \norm{U(t_{k+1}) - \temp{k+1} }_{\fh,2} \\
& \quad \le 
\br[\Big]{ \norm{ U(t_k) - \temp{k} }_{\fh,2} + (\tau_k)^2 \br[\big]{ 2\norm[]{ A^2 \temp{k} }_{\fh,2} + \norm[]{ A F(\temp{k+1}) }_{\fh,2} } } \exp(\tau_k \LL_{t_{k+1}})
\end{split}
\end{equation}
(cf.\ \cref{def:euclid_norm}).
\end{athm}

\begin{aproof}
First, note that \cref{lem:space_conv2a} ensures that for all $k \in \N_0$ with $t_{k+1} \in [0,\Ts)$ it holds that
\begin{align}\label{eq:conv1}
U(t_{k+1})
& = \exp(t_{k+1}A) U(0) + \int_0^{t_{k+1}} \exp\pr[\big]{ (t_{k+1}-s) A } F(U(s)) \dx s \nonumber \\
& = \exp\pr[\big]{ (\tau_k+t_k) A} U(0) + \int_0^{\tau_k+t_k} \exp\pr[\big]{ (\tau_k+t_k-s) A } F(U(s)) \dx s \nonumber \\
& = \exp(\tau_k A) \exp(t_k A) U(0) + \exp(\tau_k A) \int_0^{t_k} \exp\pr[\big]{ (t_k-s) A } F(U(s)) \dx s \nonumber \\
& \qquad + \int_{t_k}^{\tau_k+t_k} \exp\pr[\big]{ (\tau_k+t_k-s) A } F(U(s)) \dx s \\
& = \exp(\tau_k A) \br*{ \exp(t_k A) U(0) + \int_0^{t_k} \exp\pr[\big]{ (t_k-s) A } F(U(s)) \dx s } \nonumber \\
& \qquad + \int_{0}^{\tau_k} \exp\pr[\big]{ (\tau_k-s) A } F(U(s+t_k)) \dx s \nonumber \\
& = \exp(\tau_kA)U(t_k) + \int_{0}^{\tau_k} \exp\pr[\big]{ (\tau_k-s) A } F(U(s+t_k)) \dx s \nonumber 
\end{align}
(cf.\ \cref{def:mat_exp}).
This and \cref{lem:mon1a} assure that for $k \in \N_0$ with $t_{k+1} \in [0,\Ts)$ it holds that
\begin{align}
& U(t_{k+1}) - \temp{k+1} \nonumber \\
& = \exp(\tau_kA)U(t_k) + \int_{0}^{\tau_k} \exp\pr[\big]{ (\tau_k-s) A } F(U(s+t_k)) \dx s \nonumber \\
& \quad - \br[\big]{ (\idMatrix_N - \tau_k A)^{-1} \temp{k} + \tau_k F(\temp{k+1}) } \nonumber \\
& = \exp(\tau_kA)U(t_k) + \int_{0}^{\tau_k} \exp\pr[\big]{ (\tau_k-s) A } F(U(s+t_k)) \dx s \nonumber \\
& \quad - \br*{ \exp(\tau_kA) \temp{k} + \int_{0}^{\tau_k} \exp\pr[\big]{ (\tau_k-s) A } F(\temp{k+1}) \dx s } \\
& \quad + \br*{ \exp(\tau_kA)\temp{k} + \int_{0}^{\tau_k} \exp\pr[\big]{ (\tau_k-s) A } F(\temp{k+1}) \dx s } 
- \br[\big]{ (\idMatrix_N - \tau_k A)^{-1} \temp{k} + \tau_k F(\temp{k+1}) } \nonumber \\
& = \exp(\tau_kA) \pr[\big]{ U(t_k) - \temp{k} } + \int_{0}^{\tau_k} \exp\pr[\big]{ (\tau_k-s) A } \pr[\big]{ F(U(s+t_k)) - F(\temp{k+1}) } \dx s \nonumber \\
& \quad + \pr[\big]{ \exp(\tau_k A) - (\idMatrix_N - \tau_k A)^{-1} } \temp{k} + \int_{0}^{\tau_k} \br[\Big]{ \exp\pr[\big]{ (\tau_k-s) A } - \idMatrix_N } F(\temp{k+1}) \dx s \nonumber
\end{align}
(cf.\ \cref{def:identityMatrix}).
Combining this and the triangle inequality demonstrates that for all $k\in\N_0$ with $t_{k+1} \in [0,\Ts)$ it holds that
\begin{align}\label{eq:conv_comp}
\norm{ U(t_{k+1}) - \temp{k+1} }_{\fh,2} 
& \le \norm[\Big]{ \exp(\tau_kA) \pr[\big]{ U(t_k) - \temp{k} } }_{\fh,2} \nonumber \\
& \qquad + \norm[\big]{ \pr[\big]{ \exp(\tau_k A) - (\idMatrix_N - \tau_k A)^{-1} } \temp{k} }_{\fh,2} \\
& \qquad + \int_{0}^{\tau_k} \norm[\Big]{ \exp\pr[\big]{ (\tau_k-s) A } \pr[\big]{ F(U(s+t_k)) - F(\temp{k+1}) } }_{\fh,2} \dx s \nonumber \\
& \qquad + \int_{0}^{\tau_k} \norm[\Big]{ \br[\Big]{ \exp\pr[\big]{ (\tau_k-s) A } - \idMatrix_N } F(\temp{k+1}) }_{\fh,2} \dx s \nonumber
\end{align}
(cf.\ \cref{def:euclid_norm}).
Next, observe that \cref{lem:A_log} and \cref{lem:exp_bd} assure that for all $k\in\N_0$ with $t_{k+1} \in [0,\Ts)$ it holds that
\begin{equation}\label{eq:log_exp_norm}
\begin{split}
\norm[\Big]{ \exp(\tau_kA) \pr[\big]{ U(t_k) - \temp{k} } }_{\fh,2} 
& \le \exp(\tau_k \logNorm_{\fh,2}(A)) \norm[\big]{ U(t_k) - \temp{k} }_{\fh,2} \\
& \le \norm[\big]{ U(t_k) - \temp{k} }_{\fh,2}
\end{split}
\end{equation}
(cf.\ \cref{def:log}).
Combining \cref{eq:log_exp_norm,f_cond,eq:conv_comp} with \cref{lem:pre_converge_i2,lem:pre_converge_i2a} in \cref{lem:pre_converge} hence proves that for all $k\in\N_0$ with $t_{k+1} \in [0,\Ts)$ it holds that
\begin{align}
\norm{ U(\tau_k + t_{k}) - \temp{k+1} }_{\fh,2} 
& \le \norm{ U(t_k) - \temp{k} }_{\fh,2} + 2(\tau_k)^2 \norm[]{ A^2 \temp{k} }_{\fh,2} + (\tau_k)^2 \norm[]{ A F(\temp{k+1}) }_{\fh,2} \nonumber \\
& \qquad + \int_{0}^{\tau_k} \LL_{s+t_k} \norm{ U(s+t_k) - \temp{k+1} }_{\fh,2} \dx s \dpp
\end{align}
This, the fact that $\LL$ is a non-decreasing function, and Gronwall's inequality (cf., e.g., Henry \cite[Lemma 7.1.1]{h81}) show that for all $k\in\N_0$ such that $t_{k+1} \in [0,\Ts)$ it holds that
\begin{equation}
\begin{split}
& \norm{ U(\tau_k + t_{k}) - \temp{k+1} }_{\fh,2} \\
& \quad \le 
\br[\Big]{ \norm{ U(t_k) - \temp{k} }_{\fh,2} + (\tau_k)^2 \br[\big]{ 2\norm[]{ A^2 \temp{k} }_{\fh,2} + \norm[]{ A F(\temp{k+1}) }_{\fh,2} } } \exp(\tau_k \LL_{t_{k+1}}) \dpp
\end{split}
\end{equation}
\end{aproof}

\begin{athm}{corollary}{lem:converge2}
Assume \cref{setting3}, let $\fh\in\R^{N+1}$ satisfy $\fh = (h_0,h_1,\ldots,h_N)$, and let $\LL \colon [0,\Ts) \to [0,\infty)$ satisfy for all $t\in[0,\Ts)$ that $\LL_t = \max_{n\in\{1,2,\ldots,N\}} L_{U_n(t),\temp{\min\{k\in\N_0 \colon t \le t_k\}}}$.
Then for all $k \in \N_0$ with $t_{k+1} \in [0,\Ts)$ it holds that
\begin{equation}\label{eq:cor_converge1}
\begin{split}
& \norm{U(t_{k+1}) - \temp{k+1} }_{\fh,2} \\
& \quad \le 
\pr[\Big]{ \br[\big]{ 2\norm[]{ A^2 \temp{k} }_{\fh,2} + \norm[]{ A F(\temp{k+1}) }_{\fh,2} } \exp\pr[\big]{ {\textstyle\sum_{j=0}^k \tau_j \LL_{t_{j+1}}} } {\smallsum_{j=0}^k \tau_j } } \max_{j\in\{0,1,\ldots,k\}} \tau_j \\
& \quad \le \pr[\Big]{ \br[\big]{ 2\norm[]{ A^2 \temp{k} }_{\fh,2} + \norm[]{ A F(\temp{k+1}) }_{\fh,2} } \exp\pr[\big]{ {\textstyle\sum_{j=0}^k \tau_j \LL_{t_{j+1}} } } {\smallsum_{j=0}^k \tau_j} } \tau_0
\end{split}
\end{equation}
(cf.\ \cref{def:euclid_norm}).
\end{athm}

\begin{aproof}
Throughout this proof let $\mm_0,\mm_1,\mm_2 ,\ldots \in [0,\infty]$ satisfy for all $k \in \N_0$ that
\begin{equation}\label{eq:mm_k_fn}
\mm_k = \max_{j \in \{0,1,\ldots,k\}} \br[\big]{ 2\norm[]{ A^2 \temp{j} }_{\fh,2} + \norm[]{ A F(\temp{j+1}) }_{\fh,2} } \dpp
\end{equation}
We claim that for all $k \in \N_0$ with $t_{k+1} \in [0,\Ts)$ it holds that
\begin{equation}\label{eq:cor_converge}
\norm{U(t_{k+1}) - \temp{k+1} }_{\fh,2}
\le 
\pr[\Big]{ \mm_k \exp\pr[\big]{ {\textstyle\sum_{j=0}^k \tau_j \LL_{t_{j+1}}} } {\smallsum_{j=0}^k \tau_j} } \max_{j\in\{0,1,\ldots,k\}} \tau_j \dpp
\end{equation}
We prove \cref{eq:cor_converge} by induction.
For the base case $k=0$, note that \cref{lem:converge1} proves that
\begin{align}
& \norm{U(t_{1}) - \temp{1} }_{\fh,2} \le 
\br[\big]{ \norm{ U(t_0) - \temp{0} }_{\fh,2} + (\tau_0)^2 ( 2\norm[]{ A^2 \temp{0} }_{\fh,2} + \norm[]{ A F(\temp{1}) }_{\fh,2} ) } \exp(\tau_0 \LL_{t_1}) \nonumber \\
& \quad = (\tau_0)^2 ( 2\norm[]{ A^2 \temp{0} }_{\fh,2} + \norm[]{ A F(\temp{1}) }_{\fh,2} ) \exp(\tau_0 \LL_{t_1}) 
\le \mm_0 \tau_0 \exp(\tau_0 \LL_{t_1}) \tau_0 \dpp 
\end{align}
This establishes \cref{eq:cor_converge} in the base case $k=0$.
For the induction step $\N_0 \ni (k-1) \induct k \in \N$, let $k\in\N$ such that $t_{k+1} \in [0,\Ts)$ and assume that for every $i \in \{0,1,\ldots,k-1\}$ with $t_{i+1} \in [0,\Ts)$ it holds that
\begin{equation}
\norm{U(t_{i+1}) - \temp{i+1} }_{\fh,2}
\le 
\pr[\Big]{ \mm_i \exp\pr[\big]{ {\textstyle\sum_{j=0}^i \tau_j \LL_{t_{j+1}}} } {\smallsum_{j=0}^i \tau_j} } \max_{j\in\{0,1,\ldots,i\}} \tau_j \dpp
\end{equation}
This, 
the fact that \cref{eq:mm_k_fn} assures that for all $k\in\N_0$ it holds that $\mm_{k+1} \ge \mm_k$,
\cref{lem:mon2}, and \cref{lem:converge1} imply that
\begin{align}\label{eq:4_16}
& \norm{U(t_{k+1}) - \temp{k+1} }_{\fh,2} \nonumber \\
& \quad
\le 
\br[\big]{ \norm{ U(t_k) - \temp{k} }_{\fh,2} + (\tau_k)^2 ( 2\norm[]{ A^2 \temp{k} }_{\fh,2} + \norm[]{ A F(\temp{k+1}) }_{\fh,2} ) } \exp(\tau_k \LL_{t_{k+1}}) \nonumber \\
& \quad \le \biggl[ \pr[\Big]{ \mm_{k-1} \exp\pr[\big]{ {\textstyle\sum_{j=0}^{k-1} \tau_j \LL_{t_{j+1}} } } {\smallsum_{j=0}^{k-1}\tau_j} } \max_{j\in\{0,1,\ldots,k-1\}} \tau_j \nonumber \\
& \quad \qquad + (\tau_k)^2 \pr[\big]{ 2\norm[]{ A^2 \temp{k} }_{\fh,2} + \norm[]{ A F(\temp{k+1}) }_{\fh,2} } \biggr] \exp(\tau_k \LL_{t_{k+1}}) \\
& \quad \le \mm_k \biggl[ \pr[\Big]{ \exp\pr[\big]{ {\textstyle\sum_{j=0}^{k-1} \tau_j \LL_{t_{j+1}} } } {\smallsum_{j=0}^{k-1}\tau_j} } \max_{j\in\{0,1,\ldots,k-1\}} \tau_j + (\tau_k)^2 \biggr] \exp(\tau_k \LL_{t_{k+1}}) \nonumber \\
& \quad \le 
\pr[\Big]{ \mm_k \exp\pr[\big]{ {\textstyle\sum_{j=0}^k \tau_j \LL_{t_{j+1}} } } {\smallsum_{j=0}^k \tau_j} } \max_{j\in\{0,1,\ldots,k\}} \tau_j \dpp \nonumber 
\end{align}
Induction thus establishes \cref{eq:cor_converge}.
Next, note that \cref{eq:time_step} 
ensures that for all $k\in\N_0$ 
it holds that $\tau_k \le \tau_0$.
Combining this and \cref{eq:4_16} proves \cref{eq:cor_converge1}.
\end{aproof}

We now close \cref{sec:5_2} with some informal remarks on the bounds presented in \cref{eq:cor_converge1} in \cref{lem:converge2} above.
Essentially, the validity of the bounds in \cref{eq:cor_converge1} in \cref{lem:converge2} above (in terms of utilizing them for convergence analysis) depend upon
\begin{enumerate}[label=(\Roman*)]
\item\label{valid1} it holding for all $k \in \N_0$ that $2\norm[]{ A^2 \temp{k} }_{\fh,2} + \norm[]{ A F(\temp{k+1}) }_{\fh,2} < \infty$,
\item\label{valid2} it holding for all $k \in \N_0$ that $\smallsum_{j=0}^k \tau_j < \infty$, and
\item\label{valid3} it holding for all $k \in \N_0$ with $t_{k+1} \in [0,\Ts)$ that $\exp\pr[\big]{ {\textstyle\sum_{j=0}^k \tau_j \LL_{t_{j+1}} } } < \infty$.
\end{enumerate}
Note that \cref{valid1} is a discrete version of a regularity condition on the fully-discrete numerical solution.
We will explore this issue in more detail in \cref{th:final} below, but 
for the sake of brevity we leave major explorations of this bound for future endeavors.
Next, observe that \cref{valid2} is related to how well the numerical quenching time approximates the actual quenching time $T \in (0,\infty)$ of the solution $u$ to \cref{pde2} (cf.\ \cref{setting1}).
However, this bound is an immediate consequence of the conditions imposed in \cref{eq:time_step} and the assumptions on the nonlinear function $f \in C^1([0,1),\R)$ (cf.\ \cref{f_cond2} in \cref{setting1}).

Finally, we note that \cref{valid3} is the more subtle of the conditions above and is the point where we make the more informal remarks for the sake of clarity.
It should be clear from \cref{eq:time_step} in \cref{setting3} that for all $k \in \N_0$ it holds that $\tau_k \LL_{t_{k+1}} \approx \delta$ (this can be made more precise with improved assumptions on the temporal grid).
Thus, the boundedness of the terms in \cref{valid3} for all $k\in\N_0$ is not expected to hold, in general, due to the nature of the quenching singularity (which is measured by the function $\LL \colon [0,\Ts) \to [0,\infty)$).
However, in practice one would not allow $k \to \infty$, as round-off errors would become an issue before the validity of the bound in \cref{valid3}.

\subsection{Convergence analysis of the fully-discretized method}\label{sec:5_3}

In this subsection, we present the main results of the article. 
In \cref{th:final} below, we demonstrate that under certain reasonable continuous and discrete regularity assumptions, it holds that the implicit nonlinear operator splitting algorithm provides a numerical approximation to the solution $u$ of \cref{pde2} which has a first-order global convergence rate in space and time, prior to the point of quenching.
Moreover, prior to quenching, the nonlinear operator splitting algorithm produces approximations to the solution $u$ of \cref{pde2} which are componentise nonnegative and monotonically increasing.
The recovery of these qualitative features was an important goal of the work, herein.
Next, in \cref{cor:final} below, we apply the results from \cref{th:final} to a particular quenching-combustion PDE of interest.
In this result, it is also shown that the assumptions presented in \cref{setting1}, \cref{setting2}, \cref{setting3}, and \cref{th:final} are reasonable and easily satisfied.

\begin{athm}{theorem}{th:final}
Assume \cref{setting3}, 
let $\U \colon [0,\Ts) \allowbreak \to \R^N$ satisfy for all $t\in[0,\Ts)$, $n\in\{1,2,\ldots,N\}$ that $\U_n(t) = u(t,x_n)$, 
let $\fh \in \R^{N+1}$ satisfy $\fh = (h_0,h_1,\ldots,h_N)$,
let $\LL^{(i)} \colon [0,\Ts) \to [0,\infty)$, $i\in\{1,2\}$, satisfy for all $t\in[0,\Ts)$ that 
\begin{equation}
\LL^{(1)}_t = \max_{n\in\{1,2,\ldots,N\}} L_{u(t,x_n),U_n(t)}
\qquad \text{and} \qquad
\LL^{(2)}_t = \max_{n\in\{1,2,\ldots,N\}} L_{U_n(t),\temp{\min\{k\in\N_0 \colon t \le t_k\}}} \dc
\end{equation}
assume that $u \in C^{1,3}([0,T)\times[-a,a],\R)$, 
assume that $A \temp{0} + F(\temp{0}) - \tau_0 A F(\temp{0})$ is nonnegative,
assume that $\norm{A^2\temp{0}}_{\fh,2} < \infty$,
and assume that there exist $\kappa_0,\kappa_1,\kappa_2,\ldots \in \R$ such that for all $k \in \N_0$ it holds that $\tau_k \kappa_k \in [0,1)$ and $\max\{\norm{A F(\temp{k})}_{\fh,2},\norm{A^2 F(\temp{k})}_{\fh,2}\} \le \kappa_k \norm{A^2\temp{k}}_{\fh,2}$
(cf.\ \cref{def:wedge,def:euclid_norm}).
Then
\begin{enumerate}[label=(\roman*)]
\item\label{th:final_i1} it holds for all $k \in \N_0$, $n\in\{1,2,\ldots,N\}$ that $0 \le \temp{k}_n \le \temp{k+1}_n < 1$,
\item\label{th:final_i2} it holds for all $k\in\N_0$ that $\textstyle\sum_{j=0}^k \tau_j < \infty$,
\item\label{th:final_i3} it holds for all $k\in\N_0$ that
\begin{equation}
2\norm[\big]{ A^2 \temp{k} }_{\fh,2} + \norm[\big]{ A F(\temp{k+1}) }_{\fh,2} \le (2+\kappa_{k+1}) \exp\pr[\big]{ {\textstyle\sum_{j=0}^{k} \tau_j\kappa_j } } \norm[\big]{ A^2 \temp{0} }_{\fh,2} \dc
\end{equation}
and
\item\label{th:final_i4} it holds that there exist $\mathfrak{C}_0, \mathfrak{C}_1, \mathfrak{C}_2, \ldots \in \R$ such that for all $k \in \N_0$ with $t_{k+1} \in [0,\Ts)$ it holds that
\begin{equation}
\begin{split}
\norm{\U(t_k) - \temp{k}}_{\fh,2} 
& \le
\mathfrak{C}_k \biggl[ \exp(\Ts \LL^{(1)}_{t_k}) \max_{n\in\{1,2,\ldots,N\}} h_n \\
& \qquad  + \pr[\Big]{ (2+\kappa_{k+1}) \exp\pr[\big]{ {\textstyle\sum_{j=0}^k \tau_j ( \kappa_j + \LL^{(2)}_{t_{j+1}} ) } } {\smallsum_{j=0}^k \tau_j} } \tau_0 \biggr] \dpp
\end{split}
\end{equation}
\end{enumerate}
\end{athm}

\begin{aproof}
Throughout this proof let $\mathfrak{C}_0, \mathfrak{C}_1, \mathfrak{C}_2, \ldots \in \R$ satisfy for all $k \in \N_0$ with $t_{k+1} \in [0,\Ts)$ that
\begin{equation}\label{const_fn_bd_th}
\mathfrak{C}_k = \Ts \pr[\Bigg]{ \sqrt{2a} \left[ \max_{(t,x) \in [0,t_k] \times [-a,a]} \abs[\Big]{ \pr[\big]{\tfrac{\partial^3}{\partial x^3} u}\lrSpace(t,x) } \right] + \norm[\big]{ A^2 \temp{0} }_{\fh,2} } \dpp
\end{equation}
Note that combining \cref{lem:sol_sd_pos}, \cref{lem:fd_sol_bd}, \cref{lem:mon2}, and the assumption that $A \temp{0} + F(\temp{0}) - \tau_0 A F(\temp{0})$ is nonnegative establishes \cref{th:final_i1}.
Next, we claim that for all $a \in [0,1)$, $b \in (0,1)$ it holds that
\begin{equation}\label{eq:time_series_bd_pre}
0 \le \sum_{j=0}^\infty \frac{ab^j}{f(1-ab^j)} \le -\frac{1}{\ln b} \int_0^1 \frac{1}{f(1-s)} \dx s < \infty \dpp
\end{equation}
Observe that the assumption that $f(0) > 0$, the assumption that for all $x \in [0,1)$ it holds that $f'(x) > 0$, the fact that $(0,1] \ni s \mapsto \nicefrac{s}{f(1-s)} \in [0,\infty)$ is non-decreasing, and the fact that for all $j \in \N_0$, $s \in [j,j+1]$, $a \in [0,1)$, $b \in (0,1)$ it holds that $ab^{s} \le ab^j$ ensure that for all $j \in \N_0$, $a \in [0,1)$, $b \in (0,1)$ it holds that
\begin{equation}
\frac{ab^{j}}{f(1-ab^{j})} \le \int_j^{j+1} \frac{ab^s}{f(1-ab^s)} \dx s \dpp
\end{equation}
This assures that for all $a \in [0,1)$, $b \in (0,1)$ it holds that
\begin{equation}
\begin{split}
\sum_{j=0}^\infty \frac{ab^{j}}{f(1-ab^{j})} 
& \le \sum_{j=0}^\infty \br*{ \int_j^{j+1} \frac{ab^s}{f(1-ab^s)} \dx s }
= \int_0^{\infty} \frac{ab^s}{f(1-ab^s)} \dx s  \\
& 
= -\frac{1}{\ln b} \int_0^a \frac{1}{f(1-s)} \dx s
\le -\frac{1}{\ln b} \int_0^1 \frac{1}{f(1-s)} \dx s \dpp
\end{split}
\end{equation}
Combining this, the assumption that $f(0) > 0$, and \cref{f_cond2} proves \cref{eq:time_series_bd_pre}.
In addition, note that the fact that $[0,1) \ni s \mapsto \nicefrac{(1-s)}{f(s)} \in (-\infty,0]$ is non-increasing and \cref{lem:fd_sol_bd} demonstrate that for all $j \in \N_0$ it holds that
\begin{equation}
\begin{split}
\max_{i\in\{1,2,\ldots,N\}} \frac{1-\temp{j+1}_i  }{f(\temp{j+1}_i)}
& \le \max_{i\in\{1,2,\ldots,N\}} \frac{ (1+\delta)^{-(j+1)}(1-\temp{0}_i) }{f\pr[\big]{ 1 - (1+\delta)^{-(j+1)}(1-\temp{0}_i) }} \\
& \le \frac{ (1+\delta)^{-(j+1)}(1- {\textstyle\min_{i \in \{1,2,\ldots,N\}}}\temp{0}_i) }{f\pr[\big]{ 1 - (1+\delta)^{-(j+1)}(1-{\textstyle\min_{i \in \{1,2,\ldots,N\}}}\temp{0}_i) }} \dpp
\end{split}
\end{equation}
This, the assumption that $f(0) > 0$, the assumption that $\delta \in (0,\min\{1,\Ts\})$, and \cref{eq:time_series_bd_pre} yield that for all $k \in \N_0$ it holds that
\begin{align}
{\textstyle\sum_{j=0}^k \tau_j}
& \le {\textstyle\sum_{j=0}^\infty \tau_j}
\le \sum_{j=0}^\infty \br*{ \max_{i\in\{1,2,\ldots,N\}} \frac{\delta (1-\temp{j+1}_i ) }{f(\temp{j+1}_i)} } \nonumber \\
& \le - \frac{(1- {\textstyle\min_{i \in \{1,2,\ldots,N\}}}\temp{0}_i)}{f\pr[\big]{{\textstyle\min_{i \in \{1,2,\ldots,N\}}}\temp{0}_i)}} + \sum_{j=0}^\infty \frac{ (1+\delta)^{-j}(1- {\textstyle\min_{i \in \{1,2,\ldots,N\}}}\temp{0}_i) }{f\pr[\big]{ 1 - (1+\delta)^{-j}(1-{\textstyle\min_{i \in \{1,2,\ldots,N\}}}\temp{0}_i) }} \\
& \le - \frac{(1- {\textstyle\min_{i \in \{1,2,\ldots,N\}}}\temp{0}_i)}{f\pr[\big]{{\textstyle\min_{i \in \{1,2,\ldots,N\}}}\temp{0}_i)}} + \frac{1}{\ln(1+\delta)} \int_0^1 \frac{1}{f(1-s)} \dx s
< \infty \dpp \nonumber
\end{align}
This establishes \cref{th:final_i2}.
Next, we claim that for all $k \in \N_0$ it holds that
\begin{equation}
\label{eq:reg_claim}
\norm[\big]{ A^2 \temp{k} }_{\fh,2} \le \br*{ \prod_{j=0}^k \pr[\big]{ 1 - \tau_{\max\{j-1,0\}} \kappa_{\max\{j-1,0\}} }^{-1} } \norm[\big]{ A^2 \temp{0} }_{\fh,2} \dpp
\end{equation}
We prove \cref{eq:reg_claim} by induction on $k\in\N_0$.
Note that the base case $k=0$ holds trivially.
For the induction step $\N_0 \ni (k-1) \induct k \in \N$, let $k \in \N$ and assume that for all $i \in \{0,1,\ldots,k-1\}$ it holds that
\begin{equation}
\label{eq:reg_claim1}
\norm[\big]{ A^2 \temp{i} }_{\fh,2} \le \br*{ \prod_{j=0}^i \pr[\big]{ 1 - \tau_{\max\{j-1,0\}} \kappa_{\max\{j-1,0\}} }^{-1} } \norm[\big]{ A^2 \temp{0} }_{\fh,2} \dpp
\end{equation}
Note that \cref{lem:mon1a_i2} in \cref{lem:mon1a}, the triangle inequality, \cref{lem:A_log}, \cref{lem:exp_bd_i1} in \cref{lem:exp_bd},
and the assumption that for all $i \in \N_0$ it holds that $\norm{A^2 F(\temp{i})}_{\fh,2} \le \kappa_i \norm{A^2\temp{i}}_{\fh,2}$ guarantee that
\begin{equation}
\begin{split}
\norm[\big]{ A^2 \temp{k} }_{\fh,2}
& = \norm[\Big]{ A^2 \br[\big]{ (\idMatrix_N - \tau_{k-1} A)^{-1} \temp{k-1} + \tau_{k-1} F(\temp{k}) } }_{\fh,2} \\
& \le \norm[\big]{ (\idMatrix_N - \tau_{k-1} A)^{-1} A^2 \temp{k-1} }_{\fh,2} + \tau_{k-1} \norm[\big]{ A^2 F(\temp{k}) }_{\fh,2} \\
& \le \norm[\big]{ A^2 \temp{k-1} }_{\fh,2} + \tau_{k-1} \kappa_{k-1} \norm[\big]{ A^2 \temp{k} }_{\fh,2} 
\end{split}
\end{equation}
(cf.\ \cref{def:identityMatrix}).
Combining this, the assumption that for all $i \in \N_0$ it holds that $\tau_i \kappa_i \in [0,1)$, and \cref{eq:reg_claim1} demonstrates that
\begin{equation}
\begin{split}
\norm[]{ A^2 \temp{k} }_{\fh,2}
& \le \pr[\big]{ 1 - \tau_{k-1} \kappa_{k-1} }^{-1} \norm[]{ A^2 \temp{k-1} }_{\fh,2} \\
& \le \pr[\big]{ 1 - \tau_{k-1} \kappa_{k-1} }^{-1} \pr*{ \br*{ \prod_{j=0}^{k-2} \pr[\big]{ 1 - \tau_{\max\{j-1,0\}} \kappa_{\max\{j-1,0\}} }^{-1} } \norm[\big]{ A^2 \temp{0} }_{\fh,2} } \\
& = \br*{ \prod_{j=0}^k \pr[\big]{ 1 - \tau_{\max\{j-1,0\}} \kappa_{\max\{j-1,0\}} }^{-1} } \norm[\big]{ A^2 \temp{0} }_{\fh,2}
\dpp
\end{split}
\end{equation}
Induction hence establishes \cref{eq:reg_claim}.
In addition, observe that \cref{eq:reg_claim} and the assumption that for all $i \in \N_0$ it holds that $\tau_i \kappa_i \in [0,1)$ imply that for all $k \in \N_0$ it holds that
\begin{equation}\label{5_39}
\norm[]{ A^2 \temp{k} }_{\fh,2}
\le 
\exp\pr[\big]{ {\textstyle\sum_{j=0}^{k-1} \tau_j\kappa_j } } \norm[\big]{ A^2 \temp{0} }_{\fh,2} \dpp
\end{equation}
This and the assumption that for all $k \in \N_0$ it holds that $\norm{A F(\temp{k})}_{\fh,2} \le \kappa_k \norm{A^2\temp{k}}_{\fh,2}$
ensure that for all $k \in \N_0$ it holds that
\begin{equation}
\begin{split}
\norm[\big]{ A F(\temp{k+1}) }_{\fh,2} 
& \le
\kappa_{k+1} \norm[\big]{A^2\temp{k+1}}_{\fh,2}
\le \kappa_{k+1} \exp\pr[\big]{ {\textstyle\sum_{j=0}^{k} \tau_j\kappa_j } } \norm[\big]{ A^2 \temp{0} }_{\fh,2} \dpp
\end{split}
\end{equation}
Combining this and \cref{5_39} establishes \cref{th:final_i3}.
Next, note that the triangle inequality, \cref{lem:space_conv}, and \cref{lem:converge2} assure that for all $k \in \N_0$ with $t_{k+1} \in [0,\Ts)$ it holds that
\begin{align}\label{5_42}
\norm[\big]{\U(t_k) - \temp{k}}_{\fh,2} 
& \le \norm[\big]{\U(t_k) - U(t_k)}_{\fh,2} + \norm[\big]{U(t_k) - \temp{k}}_{\fh,2} \nonumber \\
& \le \sqrt{2a} \exp(\Ts \LL^{(1)}_{t_k}) \left[ \int_0^{t_k} \max_{n\in\{1,2,\ldots,N\}} \abs*{ \int_0^{h_n} \pr[\big]{ \tfrac{\partial^3}{\partial x^3} u}\lrSpace(w,s) \dx s } \dx w \right] \\
& \quad + \pr[\Big]{ \br[\big]{ 2\norm[]{ A^2 \temp{k} }_{\fh,2} + \norm[]{ A F(\temp{k+1}) }_{\fh,2} } \exp\pr[\big]{ {\textstyle\sum_{j=0}^k \tau_j \LL^{(2)}_{t_{j+1}} } } {\smallsum_{j=0}^k \tau_j} } \tau_0 \dpp \nonumber
\end{align}
In addition, observe that the assumption that $u \in C^{1,3}([0,T)\times[-a,a],\R)$ and Jensen's inequality show that for all $k \in \N_0$ with $t_{k+1} \in [0,\Ts)$ it holds that
\begin{equation}\label{5_43}
\begin{split}
& \int_0^{t_k} \max_{n\in\{1,2,\ldots,N\}} \abs*{ \int_0^{h_k} \pr[\big]{ \tfrac{\partial^3}{\partial x^3} u}\lrSpace(w,s) \dx s } \dx w \\
& \quad \le \int_0^{t_k} \left[ \max_{n\in\{1,2,\ldots,N\}} \int_{0}^{h_n} \max_{(t,x) \in [0,t_k] \times [-a,a]} \abs[\Big]{ \pr[\big]{\tfrac{\partial^3}{\partial x^3} u}\lrSpace(t,x) } \dx s \right] \dx w \\
& \quad = \Ts \left[ \max_{(t,x) \in [0,t_k] \times [-a,a]} \abs[\Big]{ \pr[\big]{\tfrac{\partial^3}{\partial x^3} u}\lrSpace(t,x) } \right] \max_{n\in\{1,2,\ldots,N\}} h_n \dpp
\end{split}
\end{equation}
Combining \cref{const_fn_bd_th}, \cref{5_42}, \cref{5_43}, and the fact that $\LL^{(2)} \colon [0,\Ts) \to [0,\infty)$ is non-decreasing therefore yields that
for all $k \in \N_0$ with $t_{k+1} \in [0,\Ts)$ it holds that
\begin{align}
\norm[\big]{\U(t_k) - \temp{k}}_{\fh,2} 
& \le 
\mathfrak{C}_k \biggl[ \exp(\Ts \LL^{(1)}_{t_k}) \max_{n\in\{1,2,\ldots,N\}} h_n \nonumber \\
& \qquad  + \pr[\Big]{ (2+\kappa_{k+1}) \exp\pr[\big]{ {\textstyle\sum_{j=0}^{k} \tau_j\kappa_j } } \exp\pr[\big]{ {\textstyle\sum_{j=0}^k \tau_j \LL^{(2)}_{t_{j+1}} } } {\smallsum_{j=0}^k \tau_j} } \tau_0 \biggr] \nonumber \\
& \le 
\mathfrak{C}_k \biggl[ \exp(\Ts \LL^{(1)}_{t_k}) \max_{n\in\{1,2,\ldots,N\}} h_n \\
& \qquad  + \pr[\Big]{ (2+\kappa_{k+1}) \exp\pr[\big]{ {\textstyle\sum_{j=0}^k \tau_j ( \kappa_j + \LL^{(2)}_{t_{j+1}} ) } } {\smallsum_{j=0}^k \tau_j} } \tau_0 \biggr] \dpp \nonumber
\end{align}
This establishes \cref{th:final_i4}.
\end{aproof}

\begin{athm}{remark}{rem:final}
Note that the assumption that there exist $\kappa_0,\kappa_1,\kappa_2,\ldots \in \R$ such that for all $k \in \N_0$ it holds that $\tau_k \kappa_k \in [0,1)$ and $\max\{\norm{A F(\temp{k})}_{\fh,2},\norm{A^2 F(\temp{k})}_{\fh,2}\} \le \kappa_k \norm{A^2\temp{k}}_{\fh,2}$ in \cref{th:final} above may be viewed as a type of discrete domain-invariance condition on the nonlinear function $F$ and matrix $A$. 
The fact that we utilize a sequence of constants, as opposed to a uniform one, is due to the singular nature of the problem at hand; however, a uniform bound may be introduced by studying the problem on $[0,\Ts-\varepsilon]\subseteq[0,\Ts)$, for some $\varepsilon \in (0,\min\{1,\Ts\})$.
Moreover, it is worth observing that we may replace this condition with any other discrete regularity condition which ensures that there exist $\gamma_0, \gamma_1, \gamma_2, \ldots \in \R$ such that for all $k \in \N_0$ with $t_{k+1} \in [0,\Ts)$ it holds that $\max\{\norm{A^2\temp{k}}_{\fh,2},\norm{AF(\temp{k+1})}\} \le \gamma_k$.
\end{athm}

Next, in \cref{cor:final} below, we assume that the true solution quenching time and the semidiscrete quenching time are the same (i.e., we assume that $T = \Ts$ in the language of \cref{setting1,setting2}).
In addition, as in \cref{th:final} above, we assume that there exist $\kappa_0,\kappa_1,\kappa_2,\ldots \in \R$ such that for all $k \in \N_0$ it holds that $\tau_k \kappa_k \in [0,1)$ and $\max\{\norm{A F(\temp{k})}_{\fh,2},\norm{A^2 F(\temp{k})}_{\fh,2}\} \le \kappa_k \norm{A^2\temp{k}}_{\fh,2}$.
These assumptions are only for simplicity of presentation and we leave more explicit explorations of these issues for future endeavors.

\begin{athm}{corollary}{cor:final}
Let $a,T \in (0,\infty)$, $N \in \N$,
let $u \colon [0,T) \times [-a,a] \to \R$ satisfy for all $t \in [0,T)$, $x \in [-a,a]$ that $u(t,-a) = u(t,a) = u(0,x) = 0$, $\lim_{s\to T^-} [ \sup_{y \in [-a,a]} u(s,y) ] = 1$, and
\begin{equation}
    \pr[\big]{ \tfrac{\partial}{\partial t} u}\lrSpace(t,x) = \pr[\big]{ \tfrac{\partial^2}{\partial x^2} u}\lrSpace(t,x) + \pr[\big]{ 1-u(t,x) }^{-1} \dc
\end{equation} 
let $h_0,h_1,\ldots,h_N,x_0,x_1,\ldots,\allowbreak x_{N+1} \allowbreak \in [-a,a]$ satisfy for all $n \in \{0,1,\ldots,N\}$ that $-a = x_0 < x_1 < \ldots < x_{N+1} = a$ and $x_{n+1} - x_n = h_n$, 
let $\U \colon [0,\Ts) \allowbreak \to \R^N$ satisfy for all $t\in[0,T)$, $n\in\{1,2,\ldots,N\}$ that $\U_n(t) = u(t,x_n)$,
let $F \colon \R^N \to \R^N$ satisfy for all $X \in \R^N$, $n \in \{1,2,\ldots,N\}$ that $(F(X))_n = (1-X_n)^{-1}$, 
let $A = (A_{i,j})_{i,j \in\{1,2,\ldots,N\}} \in \R^{N\times N}$ satisfy for all $i,j,n \in \{1,2,\ldots,N\}$, $k \in \{1,2,\ldots,N-1\}$ with $\abs{i-j} \in \{2,3,\ldots,N-1\}$ that $A_{i,j} = 0$,
\begin{equation}
A_{k+1,k} = \frac{2}{h_k(h_k + h_{k+1})} \dc
\quad
A_{n,n} = \frac{-2}{h_{n-1}h_n} \dc
\quad \text{and} \quad
A_{k,k+1} = \frac{2}{h_k(h_{k-1}+h_k)} \dc
\end{equation}
let $\delta \in (0,\min\{1,T\})$,
$\tau_{-1} , \tau_0,\tau_1,\tau_2,\ldots \in [0,T)$, $t_0 , t_1, t_2 , \ldots \in [0,\infty)$, $\temp{0},\temp{1},\temp{2},\ldots \in \R^N$ satisfy for all $k \in \N_0$, $n \in \{1,2,\ldots,N\}$ that $\tau_{-1} = 1$, $t_0 = 0$, $t_{k+1} = t_k + \tau_k$, $\temp{0}_n = 0$,
$\tau_k 
= \delta \br[]{ \textstyle\min_{(i,j)\in\{1,2,\ldots,N\}\times\{k,k+1\}} \pr[]{ 1 - \temp{j+1}_i }^2 }
$,
and
\begin{equation}
\temp{k+1} - \temp{k} + (\tau_k)^2 AF(\temp{k+1}) = \tau_k A \temp{k+1} + \tau_k F(\temp{k+1}) \dc
\end{equation}
and assume that there exist $\kappa_0,\kappa_1,\kappa_2,\ldots \in \R$ such that for all $k \in \N_0$ it holds that $\tau_k \kappa_k \in [0,1)$ and $\max\{\norm{A F(\temp{k})}_{\fh,2},\norm{A^2 F(\temp{k})}_{\fh,2}\} \le \kappa_k \norm{A^2\temp{k}}_{\fh,2}$
(cf.\ \cref{def:euclid_norm}).
Then
\begin{enumerate}[label=(\roman*)]
\item\label{final_cor_i1} it holds for all $t_1,t_2 \in [0,T)$, $x \in [-2,2]$ with $t_1 \le t_2$ that $0 \le u(t_1,x) \le u(t_2,x) < 1$,
\item\label{final_cor_i2} it holds for all $k \in \N_0$, $n\in\{1,2,\ldots,N\}$ that $0 \le \temp{k}_n \le \temp{k+1}_n < 1$, and
\item\label{final_cor_i3} it holds that there exist $\mathfrak{C}_0 , \mathfrak{C}_1 , \mathfrak{C}_2 , \ldots \in \R$ such that for all $k \in \N_0$ with $t_{k+1} \in [0,T)$ it holds that
$
\norm{\U(t_k) - \temp{k}}_{\fh,2} \le \mathfrak{C}_k \br[]{ \tau_0 + \textstyle \max_{n\in\{0,1,\ldots,N\}} h_n }
$.
\end{enumerate}
\end{athm}

\begin{aproof}
Throughout this proof 
let $u_0 \colon [-a,a] \to [0,1)$ satisfy for all $x \in [-a,a]$ that $u_0(x) = 0$,
let $f \colon [0,1) \to \R$ be the function which satisfies for all $x \in [0,1)$ that $f(x) = (1-x)^{-1}$,
and let $L \colon [0,1) \times [0,1) \to [0,\infty)$ satisfy for all $x,y \in [0,1)$ that $L_{x,y} = (1 - \max\{x,y\})^{-1}$.
Note that the fact that for all $x \in [0,1)$ it holds that $f(x) = (1-x)^{-1}$ assures that
for all $x,y\in[0,1)$ it holds that 
\begin{equation}\label{eq:final_cor_1}
\abs{f(x) - f(y)} 
= \abs[\big]{ (1-x)^{-1} - (1-y)^{-1} } 
\le L_{x,y} \abs{x-y} \dc
\end{equation}
\begin{equation}\label{eq:final_cor_2}
f(0) = 1 > 0 \dc
\qquad
f'(x) = (1-x)^{-2} > 0 \dc
\qquad
\lim_{w\to 1^-} f(w) = \lim_{w\to 1^-} (1-w)^{-1} = \infty \dc
\end{equation}
and
\begin{equation}\label{eq:final_cor_3}
\int_0^1 f(w) \dx w = \int_0^1 (1-w)^{-1} \dx w 
= \lim_{w \to 1^-} \pr[\big]{ -\ln\abs{1-w} } = \infty \dpp
\end{equation}
Moreover, observe that the fact that $f$ is a convex function ensures that for all $x,y\in[0,1)$, $s\in[0,1]$ it holds that
\begin{equation}\label{eq:final_cor_4}
f\pr[\big]{ sx + (1-s)y } \le sf(x) + (1-s)f(y) \dpp
\end{equation}
Combining this, \cref{eq:final_cor_1}, \cref{eq:final_cor_2}, \cref{eq:final_cor_3}, 
the fact that for all $x \in [-a,a]$ it holds that $u(0,x) = u_0(x)$,
the fact that $u_0 \in C^\infty([-a,a],[0,1))$,
and \cref{lem:quench_props} establishes \cref{final_cor_i1}.
Next, note that for all $x,y \in [0,1)$ it holds that
\begin{equation}
\frac{1-x}{f(x)} = (1-x)^2
\qquad \text{and} \qquad
\frac{1}{f'(y)} = (1-y)^2 \dpp
\end{equation}
This and the assumption that for all $k \in \N_0$ it holds that $\tau_k 
= \delta \br[]{ \textstyle\min_{(i,j)\in\{1,2,\ldots,N\}\times\{k,k+1\}} \pr[]{ 1 - \temp{j+1}_i }^2 }
$ imply that for all $k \in \N_0$ it holds that
\begin{equation}\label{eq:final_cor_5}
\begin{split}
\frac{\tau_k}{\delta}
& = \min_{i \in \{1,2,\ldots,N\}} \pr*{ \min\cu*{ \pr[\big]{ 1 - \temp{k+1}_i }^2 , \pr[\big]{ 1 - \temp{k}_i }^2 } } \\
& = \min_{i \in \{1,2,\ldots,N\}} \pr*{ \min\cu*{ \frac{1-\temp{k+1}_i}{f(\temp{k+1}_i)} , \frac{1}{f'(\temp{k}_i)} } } \dpp
\end{split}
\end{equation}
In addition, observe that the assumption that for all $x \in [-a,a]$ it holds that $u(0,x) = 0$ yields that
\begin{equation}\label{eq:final_cor_6}
A\temp{0} + F(\temp{0}) - \tau_0AF(\temp{0}) = F(\temp{0}) > 0
\end{equation}
and
\begin{equation}\label{eq:final_cor_7}
\norm[\big]{A^2 \temp{0}}_{\fh,2} = 0 < \infty
\end{equation}
(cf.\ \cref{def:wedge}).
Moreover, note that the fact that $u_0 \in C^\infty([-a,a],[0,1))$ and $f \in C^\infty([0,1),\R)$ guarantee that $u \in C^{1,3}([0,T)\times[-a,a],\R)$.
Combining this, \cref{eq:final_cor_1}, \cref{eq:final_cor_2}, \cref{eq:final_cor_3}, \cref{eq:final_cor_4}, \cref{eq:final_cor_5}, \cref{eq:final_cor_6}, \cref{eq:final_cor_7}, 
and \cref{th:final}
(applied with $a \with a$, $\Ts \with T$, $u \with u$, $u_0 \with u_0$, $\U \with \U$, $L \with L$, $F \with F$, $\fh \with \fh$, $N \with N$, $A \with A$, $(\temp{k})_{k\in\N_0} \with (\temp{k})_{k\in\N_0}$, $\delta \with \delta$, $(\tau_k)_{k\in\N_0} \with (\tau_k)_{k\in\N_0}$, $(\kappa_k)_{k\in\N_0} \with (\kappa_k)_{k\in\N_0}$
in the notation of \cref{th:final})
establishes \cref{final_cor_i2,final_cor_i3}.
\end{aproof}


\section*{Acknowledgments}

The first author gratefully acknowledges funding by the National Science Foundation (NSF 1903450).

\bibliographystyle{acm}
\bibliography{Dissertation_BibTex}

\end{document}